\documentclass[12pt,a4paper]{amsart}
\usepackage[utf8]{inputenc}
\usepackage[english]{babel}
\usepackage{amsmath}
\usepackage{amsfonts}
\usepackage{amssymb}
\usepackage{graphicx}

\usepackage{tikz}
\usepackage{tikz-cd}
\usetikzlibrary{calc}
\usetikzlibrary{intersections}
\usetikzlibrary[decorations.pathmorphing,mindmap]
\usetikzlibrary{decorations.markings}
\usetikzlibrary{decorations.pathreplacing}
\usetikzlibrary{arrows.meta, decorations, positioning}

\newcommand{\QQ}{\mathbf{Q}}
\newcommand{\ZZ}{\mathbf{Z}}
\newcommand{\RR}{\mathbf R}
\newcommand{\CC}{\mathbf C}

\newcommand{\pgl}{\mathrm{PGL}_2 (\ZZ)}

\newcommand{\TT}{\begin{pmatrix}
		1 & 1\\ 0 & 1
	\end{pmatrix}
}
\newcommand{\UU}{\begin{pmatrix}
		0 & 1\\ 1 & 0
	\end{pmatrix}
}

\newtheorem{theorem}{Theorem}[section] % 1st argument is your name for it
\newtheorem{proposition}[theorem]{Proposition}
\newtheorem*{definition}{Definition}
 \newtheorem*{remark}{Remark}
\title[Thurston's Arithmetic Legacy]{Thurston's Legacy at the Interface of Surface Dynamics, Train Tracks, and Arithmetic}
\author{Ahmad Rafiqi}
\address{Department of Mathematics and Statistics\\
	American University of Sharjah\\
	Sharjah, United Arab Emirates}
\email{arafiqi@aus.edu}

\author{Ayberk Zeytin}
\address{Department of Mathematics\\
	Galatasaray University\\
	Ortaköy, Istanbul, Turkey}
\email{azeytin@gsu.edu.tr}

\subjclass[2020]{Primary 37E30, 57K20; Secondary 37B10, 30F60, 57K32, 15B48}

\keywords{multi-variate continued fractions, non-negative matrices,  pseudo-Anosov homeomorphisms, stretch factors, Teichm\"uller polynomial, Thurston norm, fibered cones, bi-Perron units, train tracks, post-critically finite maps, mapping tori.}

\thanks{The second author would like to thank IHES for their hospitality. The second author was supported by Galatasaray Research Grant FBA-2024-1236. We are thankful to the editor for their meticulous work.}

\dedicatory{In memory of William P. Thurston, with lasting admiration.}

\begin{document}
	
%	\section{Introduction}
%	The theory of non-negative matrices (a matrix is said to be non-negative (resp. positive) if all of its entries are non-negative(resp. positive).) constitutes a foundational framework that significantly influences a diverse array of disciplines, many of which appear to be conceptually unrelated, including dynamical systems, economics,etc.  The Perron–Frobenius theorem is widely recognized as one of the foundational results that plays a pivotal role in the development of this theory. Interested reader is encouraged to consult \cite{hawkins/origins/of/Perron/Frobenius/theorem} for a historical overview of the development of the result. There are other proofs of this result, see for instance \cite{wielandt/unzerlegbare/nicht/negative/matrizen, cairns/perrons/theorem}. The theory was further enriched by Lind’s proof of a converse result, alongside various generalizations—see, for instance, \cite{karlin/positive/operators}. 
%	

	\begin{abstract}		
		In this chapter, we discuss the interaction between low-dimensional dynamics, symbolic dynamics, and geometric topology, focusing on the algebraic characterization and asymptotic behaviour of stretch factors or equivalently pseudo-Anosov dilatations. We first review the algebraic foundations linking Perron--Frobenius theory, multi-variable Euclidean algorithms, and Lind’s characterization of topological entropy for shifts of finite type, emphasizing recent constraints on the Perron--Frobenius degree of expansion constants. We then relate these dynamical invariants to the combinatorial world, examining train tracks for surfaces and automorphisms of free groups. Finally, we briefly review the relevance of fibered 3-manifolds, detailing how the Thurston norm, Fried’s concavity theorem over fibered cones, and McMullen’s Teichmüller polynomial come into play. 
	\end{abstract}
	\maketitle
	\section{Introduction}

	William P. Thurston's work altered, in a decisive and permanent way, the language of low-dimensional topology and dynamical systems.  After his passing, several accounts of his mathematical legacy have appeared, underlining the breadth and profound impact of his contributions: the geometrization of $3$-manifolds, the topology and dynamics of surface homeomorphisms, measured foliations, train tracks, laminations, the Thurston norm, and the geometric structures attached to manifolds and groups; see, for example, \cite{FathiLaudenbachPoenaru2012ThurstonWork,Thurston1988GeometryDynamics,thurston1986norm,Papadopoulos2024FootnotesThurstonNotes,OhshikaPapadopoulos2020GlimpseThurston}.  The present chapter has a narrower aim.  It is not a general memorial account of Thurston's work, nor a survey of all the domains transformed by his ideas.  Its purpose is to isolate and organize one particular line of research through his contributions: the passage from topology and dynamics to finite combinatorics, and then from finite combinatorics to algebraic restrictions on spectral radii, entropies, and stretch factors.
	
	The objects discussed in this work occur under different names in several adjacent topics.  In one-dimensional dynamics one studies piecewise-monotone interval maps, Markov partitions, transition matrices, and growth rates of lap numbers.  In surface dynamics the corresponding language is that of mapping classes, invariant measured foliations, measured laminations, train tracks, and pseudo-Anosov dilatations.  In the study of free-group automorphisms one meets graph maps and transition matrices.  In $3$-manifold topology similar numerical invariants are organized by fibered faces, Fried's entropy function, and McMullen's Teichmüller polynomial.  Finally, from the algebraic point of view, the relevant numbers are Perron numbers, weak Perron numbers, bi-Perron numbers, and their degrees.  A dictionary between these languages is not straightforward.  We believe a coherent synthesis is therefore useful: it makes explicit which statements are genuinely topological, which are consequences of Perron--Frobenius theory, and which are arithmetic restrictions on algebraic integers.
	
	The central point of this chapter is: part of Thurston's contribution to entropy and stretch-factor problems has the following structure:  A geometric or topological dynamical system is first encoded by a finite combinatorial object: a Markov partition for an interval map, a train track for a pseudo-Anosov mapping class, a graph map representing an outer automorphism, or a fibered face together with its associated polynomial data.  The resulting incidence or transition matrix is non-negative and integral.  Its Perron--Frobenius eigenvalue then records the exponential growth rate of the system: topological entropy in dimension one, pseudo-Anosov dilatation on surfaces, or the stretch factor of a monodromy in a fibered $3$-manifold.  Conversely, the topology of the original system imposes restrictions on which Perron--Frobenius eigenvalues may occur.  Thus the survey is organized around the following two-way passages:
	\begin{center}
		\[
		\substack{\text{topology}\\ \text{dynamics}}
		\quad \longleftrightarrow \quad
		\substack{\text{finite}\\ \text{combinatorics}}
		\quad \longleftrightarrow \quad
		\substack{\text{algebraic}\\ \text{numbers}}.
		\]
	\end{center}
	
	The chapter begins with Perron--Frobenius theory and the algebraic properties of Perron numbers because these provide the common spectral language for the rest of the chapter.  It then treats Thurston's realization theorem for entropies of post-critically finite interval maps and its relation to Lind's theorem on Perron numbers.  The next part concerns surface homeomorphisms: the Nielsen--Thurston classification, pseudo-Anosov maps, measured foliations, train tracks, incidence matrices, and the algebraic restrictions satisfied by dilatations.  We then discuss extensions to graph maps and free-group automorphisms, least stretch factors of mapping classes, and the organization of stretch factors through fibered $3$-manifolds, the Thurston norm, Fried's theorem, and McMullen's Teichmüller polynomial. In this chapter, we use the modern theory of mapping class groups, Teichmüller dynamics, or $\mathrm{Out}(F_n)$ only to the extent needed.  Nor do we attempt an exhaustive bibliography of all recent developments on minimal dilatations; instead, we emphasize those results that clarify the relation between topology, Perron--Frobenius theory, and algebraic number theory.
	
	\section{Perron--Frobenius theory}
	
	This section supplies the common linear-algebraic language used throughout the rest of the chapter. In each of the dynamical settings considered below, a continuous, topological or algebraic object is replaced by a finite combinatorial model: a Markov partition for an interval map, a train track for a pseudo-Anosov homeomorphism, a graph train track representative of an outer automorphism of a free group, or the polynomial data associated with a fibered face of the Thurston norm ball. In all these cases one obtains a non-negative integral matrix, and the relevant dynamical invariant is its spectral radius. Thus Perron--Frobenius theory explains why topological entropy, stretch factors, dilatations, and expansion constants are governed by a distinguished positive eigenvalue. It also gives the first arithmetic restriction on these numbers: they are Perron or weak Perron algebraic integers. The later sections repeatedly use this  principle.
%	, first for post-critically finite interval maps, then for train-track representatives of surface homeomorphisms and free-group automorphisms, and finally in the discussion of stretch factors arising from fibered 3-manifolds.
	
	We now review the Perron--Frobenius theory needed for this purpose. We begin with Perron's theorem, which establishes that every positive integral matrix possesses a unique dominant eigenvalue \(\lambda\) associated with a one-dimensional eigenspace. We then discuss Frobenius's extension of this result to non-negative, irreducible and aperiodic matrices. Finally, we examine the converse direction, as proved by Lind.
	
	\subsection{The ``if'' part.}
	Given any real number $x$, the \emph{regular continued fraction}\index{continued fraction!regular} of $x$ is described as posing $x_{0} = x$, $a_{0} = \lfloor x \rfloor$ (i.e. the integer part of $x$)  and $r_{0} = x_{0} - a_{0}$ and repeating this procedure for $x_{n+1} = \frac{1}{r_{n}}$ instead of $x_{0}$, that is declaring $a_{n+1} = \lfloor x_{n+1} \rfloor$, $r_{n+1} = x_{n+1} - a_{n+1}$. Namely,
	\begin{align*}
		x & = \lfloor x \rfloor + (x - \lfloor x \rfloor)\\
		& = a_{0} + r_{0} \\
		&= a_{0} + \frac{1}{\frac{1}{r_{0}}} \\ 
		& = a_{0} + \frac{1}{x_{1}} \\
		&=  a_{0} + \frac{1}{a_{1} + r_{1}} \\
		& = a_{0} + \frac{1}{a_{1} + \frac{1}{\frac{1}{r_{1}}}} \\
		& = a_{0} + \frac{1}{a_{1} + \frac{1}{\ddots}} = [a_{0};a_{1},\ldots].
	\end{align*}
	In such an expansion, $a_{0} \in \ZZ $ and $a_{k} \in \ZZ_{+}$ for all $k = 1,2,\ldots$. The algorithm terminates if and only if $x \in \QQ$. This algorithm can be expressed using the elements 
	\[ U = \UU \quad  \mbox{and} \quad T = \TT\]
	of the group $\pgl$ considered together with its action on $\RR\cup\{\infty\}$ by M\"{o}bius transformations. Indeed, for $x \in \RR\cup \{\infty\}$, we have:
	\[ U \cdot x = \frac{1}{x} \quad \mbox{ and } \quad T^{n} \cdot x = x+n \mbox{, for }n \in \ZZ. \]
	Therefore, 
	\begin{align*}
		x &= T^{a_{0}} (r_{0}) \\
		& = T^{a_{0}} U(x_{1}) \\
		&= T^{a_{0}} U T^{a_{1}} U(x_{2}) \\
		& = \ldots 
	\end{align*}
	Equivalently, this regular continued fraction can also be interpreted as the Euclidean algorithm applied to the pair $(x,1)$. As such, given any two integers (or real numbers in the case of continued fractions), say $x_{0}$ and $x_{1}$, with $x_{0}>x_{1}$, the usual Euclidean algorithm, that is:
	\begin{align*}
		x_{0} &= x_{1} \cdot a_{0} + r_{0} \\
		x_{1} &= r_{0} \cdot a_{1} + r_{1} \\
		\vdots\quad &= \quad \vdots
	\end{align*}
	can be encoded in terms of matrices as:
	\begin{align*}
		\begin{pmatrix}
			x_{0} \\ x_{1}
		\end{pmatrix} & = T^{a_{0}} \begin{pmatrix}
		r_{0} \\ x_{1}
		\end{pmatrix}\\
								& = T^{a_{0}} U \begin{pmatrix}
									x_{1} \\ r_{0}
								\end{pmatrix} \\
								& = T^{a_{0}} U T^{a_{1}} U \begin{pmatrix}
									x_{2} \\ r_{1}
								\end{pmatrix} \\
								 &= \ldots
	\end{align*}
	In particular, this process terminates whenever $x_{0}$ and $x_{1}$ are commensurable in the sense of Euclid, \cite[Definition~1]{euclidBookX}. This interpretation leads us to express the task of finding the greatest common divisor of $(n+1)$ positive integers (or real numbers), say $x_{0},x_{1},\ldots,x_{n}$ with $x_{i}>x_{n}$ for $i = 0,1,\ldots,n-1$, using matrices. Namely, one writes the following system of equations: 
	\begin{align*}
		x_{0} &= x_{n} a_{0,0}  + r_{0,0} \\
		x_{1} &=  x_{n} a_{1,0}  + r_{1,0} \\
		\vdots &= \vdots \\
		x_{n-1} &= x_{n} a_{n-1,0}  + r_{n-1,0} \\
		x_{n} &= x_{n} \cdot 1 + 0
	\end{align*}
	in matrix form as: 
	\begin{align*}
		\begin{pmatrix}
			x_{0}\\x_{1}\\ \vdots \\x_{n-1} \\ x_{n}
		\end{pmatrix}
		& = 
		\begin{pmatrix}
			1 & 0 & 0 &  \cdots & a_{0,0} \\
			0 & 1 & 0 &  \cdots & a_{1,0} \\
			\vdots & 0 & \ddots & \vdots \\
			0 & 0 & \ldots & 1 & a_{n-1,0} \\
			0&0&0&\cdots&1
		\end{pmatrix}
		\begin{pmatrix}
			r_{0,0} \\ r_{1,0}\\ \vdots \\ r_{n-1,0}\\x_{n}
		\end{pmatrix} \\
		&=  T^{(1)}		\begin{pmatrix}
			r_{0,0} \\ r_{1,0}\\ \vdots \\ r_{n-1,0}\\x_{n}
		\end{pmatrix}
	\end{align*}
	Set $x^{(1)}_{n} = \min_{i = 0,\ldots,n-1} \{r_{i,0} \} $ and use a permutation matrix, say $U^{(1)}$, so that the last entry of the matrix 
	$$U^{(1)} \begin{pmatrix}
		r_{0,0} & r_{1,0} & \ldots & x_{n}
	\end{pmatrix}^{\top}$$
	is $x^{(1)}_{n}$. This sets the stage for an iteration of the above procedure and one has:
	\begin{align*}
				\begin{pmatrix}
			x_{0}\\x_{1}\\ \vdots  \\ x_{n}
		\end{pmatrix}
		& = T^{(1)} U^{(1)} \begin{pmatrix}
			r_{0,1} \\ r_{1,1}\\ \vdots \\x^{(1)}_{n}
		\end{pmatrix} \\
		&= T^{(1)} U^{(1)} T^{(2)}U^{(2)} \begin{pmatrix}
			r_{0,2} \\ r_{1,2}\\ \vdots \\x^{(2)}_{n}
		\end{pmatrix}.
	\end{align*}
	In the case where the $x_{i}$s are integers (or pairwise commensurable real numbers) the above algorithm terminates. The output is then a sequence of non-negative matrices: $T^{(1)}U^{(1)}\cdots T^{(k)}U^{(k)}$.% \footnote{Perron places the smallest of the numbers $x_{0},\ldots,x_{n}$ as the very last entry. This leads to slightly different matrices.} 
	This is referred to as the \emph{multi-variate continued fraction expansion}\index{continued fraction!multi-variate} of the numbers $x_{0}, \ldots, x_{n}$. By ${\bf x}$ we denote the vector $\begin{pmatrix}
		x_{0} & x_{1} & \ldots & x_{n}
	\end{pmatrix}^{\top}$.
	
	An immediate observation is that the product of matrices that arise from the above described multi-variate continued fraction algorithm is non-negative. The natural questions of convergence and periodicity therefore arise: under what hypotheses on the input vector does the associated sequence of iterates converge? Likewise, when the numbers $x_{0}, \ldots, x_{n}$ are arithmetic (for instance, lie in a number field), does the algorithm exhibit eventual periodicity, and if so, how is the period reflected in the induced matrix cocycle? It was Perron who studied, \cite{perron/theorie/jacobischen/kettenbruchalgorithmus}, convergence questions, as well as periodicity questions (in the search for an analogue of Lagrange’s theorem).\footnote{Recall that the Lagrange theorem establishes the equivalence between quadratic irrationals and periodic regular continued fraction expansions.}

	The following elementary spectral calculation shows why Perron's study of multi-variate continued fractions naturally led to the distinguished eigenvalue guaranteed by Perron--Frobenius theory (Theorems~\ref{thm:perron/root} and~\ref{thm:frobenius/theorem/for/irreducible/matrices}). Specifically, for a periodic multi-variate continued fraction expansion, the convergence of the continued fraction algorithm is equivalent to the projective convergence of the iterates of the associated matrix. In fact, Perron investigated purely periodic multi-variate continued fraction expansions by studying the linear dynamical system given by the iterates $A^{m}{\bf x}$, where $A$ encodes one period of the corresponding continued fraction, namely
	\[ A = T^{(1)}U^{(1)}\cdots T^{(k)}U^{(k)}.\]
	He observed that the convergence properties of the associated continued fraction are governed by the spectral data of $A$. More explicitly, for clarity, we assume that $A$ is diagonalizable, noting that the Jordan form yields the same asymptotic result. This allows us to express the input vector in the basis of eigenvectors $v_{0}, \ldots, v_{n}$:
	\[  {\bf x} = \begin{pmatrix}
		x_{0} \\ \vdots \\ x_{n}
	\end{pmatrix}  =  \sum_{i = 0}^{n} c_{i} v_{i}. \]
	Applying the linear operator $A^{m}$ yields $A^{m} {\bf x} = \sum_{i = 0}^{n} c_{i} \lambda_i^{m} v_{i},$
	where $\lambda_{i}$ is the eigenvalue corresponding to the eigenvector $v_{i}$. We see that for any fixed eigenvector $v_{j}$, we have:
	\begin{align}
		A^m{\bf x} = \lambda_j^m \left( c_jv_j +  \sum_{\substack{0\leq i\leq n\\ i\neq j}} c_i\left(\frac{\lambda_i}{\lambda_j}\right)^m v_i \right).
		\label{eq:geometric/interpretation}
	\end{align}
	This means that,  if \(A\) has a unique eigenvalue of maximal modulus, the projective direction of \(A^m{\bf x}\) converges to the corresponding eigenline. Perron's theorem gives precisely this conclusion for positive matrices, together with the positivity of the limiting eigenvector. We will refer to this result as \emph{Perron's theorem}\index{theorem!Perron}:
	\begin{theorem}[{\cite[\S~5]{perron/theorie/jacobischen/kettenbruchalgorithmus}}]
		Let $A$ be a positive square matrix. Then $A$ has a positive real eigenvalue. Its largest eigenvalue $\lambda_0>0$ is simple and satisfies $|\lambda_0|>|\lambda|$ for every other eigenvalue $\lambda$ of $A$. Moreover, all $n$-row minors of the matrix $\lambda_0 I-A$ are positive.
		\label{thm:perron/root}
	\end{theorem}
	In view of (\ref{eq:geometric/interpretation}), we may deduce that the orbit $\{A^{k} {\bf x}\}$ is pulled toward the so-called \emph{Perron eigenvector}\index{eigenvector!Perron}, ${\bf x}_{0}$, with asymptotic growth rate being $\lambda_{0}$, the dominant eigenvalue (or \emph{Perron root}\index{root!Perron}). 
	Notice moreover that the last assertion implies that an eigenvector associated with the Perron root $\lambda_{0}$ can be chosen with strictly positive entries. It is then easy to see that if $A$ is, more generally, a non-negative matrix with some power being positive, that is, $A^{\nu}$ is a positive matrix, then the conclusions of the above theorem are still true, as observed by Perron. Matrices satisfying the latter condition are called \emph{aperiodic}\index{matrix!aperiodic}.
	
	Frobenius realized that, for a general non-negative matrix $A$, it can be subtle to identify those eigenvalues whose eigenvectors may be chosen non-negative. He observed that the situation becomes more transparent when $A$ is similar via a permutation matrix (i.e.\ \emph{permutationally similar}\index{matrix!permutationally similar}) to a block upper triangular matrix with square diagonal blocks. In particular, if $A$ is (permutationally) similar to
	\[ 	\begin{pmatrix}
			B_1 & * & \cdots & *\\
			0 & B_2 & \cdots & *\\
			\vdots & \vdots & \ddots & \vdots\\
			0 & 0 & \cdots & B_r
		\end{pmatrix},\]	
	 then the spectrum of $A$ is the union of the spectra of the diagonal blocks $B_i$, and hence $\rho(A)=\max_{1\le i\le r}\rho(B_i)$, where $\rho(M)$ denotes the spectral radius of a matrix $M$. When the blocks are positive (or, more generally, irreducible), each $B_i$ has a distinguished Perron root that dominates the remaining eigenvalues of $B_i$ in modulus, and the spectral radius of $A$ is realized by the largest of these Perron roots. This motivates the following dichotomy: a non-negative matrix is called \emph{reducible}\index{matrix!reducible} if it is permutationally similar to a block upper triangular matrix with at least two diagonal blocks, and \emph{irreducible}\index{matrix!irreducible} otherwise. Frobenius further showed that a non-zero matrix $A$ is reducible if and only if some minor of $\lambda_0 I-A$ vanishes, where $\lambda_0=\rho(A)$ is the Perron root. The irreducible case is governed by the next theorem.
	\begin{theorem}[Frobenius Theorem, \cite{frobenius/matrizen/aus/nicht/negativen/elementen}]
		Let $A$ be a non-zero irreducible matrix. Then the largest positive eigenvalue $\lambda_0$ is simple and $|\lambda_0|\geq|\lambda|$ for every other eigenvalue $\lambda$ of $A$. Moreover, all $n$-row minors of the matrix $\lambda_0 I-A$ are positive.
		\label{thm:frobenius/theorem/for/irreducible/matrices}
	\end{theorem}
	
	Notice that, in contrast to Perron's original statement, the norm of the leading eigenvalue may coincide with norms of other eigenvalues. However, Frobenius demonstrated that such equality cannot occur if $A^{\nu}$ is positive for some $\nu\in \ZZ_{+}$, providing another proof of Perron's result.  In his work to understand the relationship between irreducibility and aperiodicity, Frobenius defined a matrix $A$ to be \emph{primitive}\index{matrix!primitive} whenever the largest positive eigenvalue \( \lambda_0 \) is simple and $|\lambda_0|>|\lambda|$ for all other eigenvalues $\lambda$ of $A$. With this terminology we have:
	\begin{theorem}[{\cite{frobenius/matrizen/aus/nicht/negativen/elementen}}]
		For an irreducible matrix $A$, primitivity is equivalent to aperiodicity.
	\end{theorem}
	This result establishes a direct correspondence between these two structural properties of non-negative matrices.
	
	\subsection{The ``only if'' part.}\label{sec:dominant/eigenvalue}
	
	An \emph{algebraic integer}\index{algebraic integer} is a complex number that is a root of a monic polynomial in $\ZZ[x]$. If $\alpha\in\CC$ is an algebraic integer, then there exists a unique monic irreducible polynomial $p_{\alpha}(x)\in\ZZ[x]$ admitting $\alpha$ as a root; this polynomial is called the \emph{minimal polynomial}\index{polynomial!minimal} of $\alpha$. The other roots of $p_{\alpha}(x)$ are called the \emph{Galois conjugates} of $\alpha$. The \emph{degree} of $\alpha$, denoted $\deg(\alpha)$, is defined to be the degree of the polynomial $p_{\alpha}$,  $\deg(p_{\alpha})$.
	
	Writing
	$p_{\alpha}(x)=a_{0}+a_{1}x+\cdots+a_{d-1}x^{d-1}+x^{d}$, one defines the \emph{companion matrix}\index{matrix!companion} of $\alpha$ (or the polynomial $p_{\alpha}$) to be
	\[ B_{\alpha}=\begin{pmatrix}
		0&0&\cdots&0&-a_{0} \\
		1&0&\cdots&0&-a_{1} \\
		0&1&\cdots&0&-a_{2} \\
		\vdots&\vdots&\ddots&\vdots&\vdots \\
		0&0&\cdots&1&-a_{d-1}
	\end{pmatrix}. \]	
	The characteristic polynomial of $B_{\alpha}$ is the minimal polynomial of $\alpha$. Therefore, $\alpha$ itself is an eigenvalue of $B_{\alpha}$. 

	Given an integer matrix $A$, its characteristic polynomial is monic with coefficients in $\ZZ$, and hence every eigenvalue of $A$ is an algebraic integer. 
	\begin{definition}
		A real algebraic integer \(\lambda>1\) is called a \emph{Perron number}\index{number!Perron} if every Galois conjugate \(\lambda'\neq \lambda\) satisfies
		\[
		|\lambda'|<\lambda .
		\]
	\end{definition}
	In fact, if $A$ satisfies the hypotheses of Perron’s theorem, then, by Theorem~\ref{thm:perron/root}, the Perron root is real and strictly positive; in this case $\lambda_{0}$ is a Perron number\index{number!Perron}. More generally, a real algebraic integer $\beta\geq 1$ is called a \emph{weak Perron number}\index{number!weak Perron} if for every Galois conjugate $\beta'$ of $\beta$ one has $|\beta'|\leq \beta$; equivalently, $\beta$ attains the spectral radius among its conjugates, possibly with multiplicity, see Theorem~\ref{thm:frobenius/theorem/for/irreducible/matrices}. In this context, Frobenius’ theorem extends Perron’s picture to irreducible non-negative integer matrices: while $\rho(A)$ may have other eigenvalues on the circle $|\lambda|=\rho(A)$ (reflecting periodicity), it is always a weak Perron number, and if $A$ is primitive (i.e.\ irreducible and aperiodic) then $\rho(A)$ is in fact a Perron number and its associated eigenvector can be chosen with strictly positive entries. Notice further that every Perron number is a weak Perron number, but not conversely.
	
	The natural question to ask is the converse, which is proven by Lind:
	\begin{theorem}[{\cite[Theorem~1]{lind/entropies/of/markov/shifts}}]
		Given any Perron number $\lambda$,  there exists an aperiodic non-negative integral matrix which admits $\lambda$ as its dominant eigenvalue. 
		\label{thm:lind}
	\end{theorem}
	
	The proof given by Lind relies on the construction of a cone around the positive direction of an eigenvector, say $x_{\lambda}$, corresponding to the eigenvalue $\lambda$ of the companion matrix of the minimal polynomial of $\lambda$. Indeed, one may find points $z_{1},\ldots,z_{n}$ which have positive $\lambda$-coordinate whose images under $ B_{\lambda}$ can be expressed as a linear combination of $z_{1},\ldots,z_{n}$ where every coefficient is a positive integer. An irreducible component of the matrix formed by these coefficients satisfies the requirements. Another proof of Theorem~\ref{thm:lind} due to Thurston will be given in Section~\ref{sec:thurston/lind/theorem}.
	
	\section{Entropies of interval maps}
	\label{sec:interval-maps}
	
	In \cite{thurston/entropy/in/dimension/one}, Thurston used piecewise-monotone interval maps as a one-dimensional model for the train-track constructions appearing later. The goal of this section is to explain how the entropy of such maps is encoded by finite combinatorial data. More precisely, after passing to a suitable Markov or post-critical partition, the dynamics is represented by a non-negative integral matrix, and the entropy is the logarithm of its spectral radius; see Sections~\ref{subsec:pcf-interval-maps} and \ref{subsec:incidence-matrices-interval-maps}.
	
	Thus Thurston's theorem has two parts: an obstruction and a realization statement. If a post-critically finite interval map has entropy \(h\), then \(\exp(h)\) is a weak Perron number; conversely, every weak Perron number arises as \(\exp(h_{\mathrm{top}})\) for some post-critically finite interval map, even for a polynomial map; see \cite[Theorem~1]{thurston/entropy/in/dimension/one} and Section~\ref{subsec:realization-powers-perron}.
	
	\subsection{Post-critically finite interval maps}
	\label{subsec:pcf-interval-maps}
	By $I$, denote a closed interval $[a,b]$ and let $f: I \to I$ be a continuous map. Say $f$ is \emph{piecewise-monotone}\index{map!piecewise-monotone} if $I$ admits a partition into sub-intervals on which $f$ is monotone. The smallest number of such intervals is called the \emph{lap number}\index{number!lap} of $f$ and each interval that appears in such a partition of $I$ is called a \emph{lap}. By $f^{n}$ we denote the composition of $f$ with itself $n$ times ($n \geq 0$, with $f^{0}$ being the identity map), and by $c_{n}$ we denote the lap number of $f^{n}$. 
	The \emph{topological entropy}\index{entropy!topological}, denoted by $h(f)$, of such a map is a fundamental measure of the asymptotic complexity of the orbit structure of the dynamical system associated with $ f $. One method to compute $h(f)$ is the context of the following result:
	\begin{theorem}[{\cite{misiurewicz/entropy/of/piecewise/linear/maps}}]
		The entropy of $f : I \to I $, denoted $h(f)$, is $\lim_{n \to \infty} \frac{\log(c_{n})}{n}$.
	\end{theorem}
	Two self-maps \(f\) and \(F\) of an interval \(I\) are \emph{semi-conjugate}\index{map!semi-conjugate} if there is a continuous onto map \(\phi:I\to I\) such that \(\phi\circ f=F\circ \phi\). If, moreover, \(\phi\) is a homeomorphism, then \(f\) and \(F\) are \emph{topologically conjugate}\index{map!topologically conjugate}. If $f$ and $F$ are topologically conjugate, then $ h(f) = h(F)$. Given $f$, its entropy $h(f)$ remains invariant under topological conjugacy. This is the motivation for one to find a function $F$, which has the same entropy as $f$ and whose entropy is easier to compute. 
	\begin{theorem}[{\cite{milnor/thurston/iterated/maps/of/the/interval}}]
		Let $f: I \to I$ be a continuous, piecewise-monotone map with a finite lap number and positive topological entropy $h(f) > 0$. Then $f$ is topologically semiconjugate to a piecewise linear map $F: I \to I$ of constant absolute slope, such that 
		\[
		|F'(x)| = \lambda = e^{h(f)} > 1
		\]
		wherever the derivative is defined. The map $F$ is called a \emph{uniform $\lambda$-expander} \index{uniform expander}, and its slope directly recovers the entropy via $h(f) = h(F) = \log \lambda$.
		\label{thm:existence/of/uniform/expander}
	\end{theorem}
	
	The proof of this result is obtained by first observing that any continuous, piecewise-monotone self-map of $I$ with a finite lap number can be associated with the symbolic dynamics defined by the sequence of laps visited under iteration. The growth rate of the number of laps of $f^n$ is used to determine the entropy, and it is noted that this combinatorial growth is independent of the precise geometric form of the laps. A monotone, continuous, surjective factor map $\phi$ is then constructed from the original interval onto another interval, in which points with identical symbolic itineraries (the itinerary records the successive laps visited by the forward orbit of $x$ under $f$, giving a symbolic representation of its dynamics) are collapsed. In this way, the ``essential'' combinatorial structure of the dynamics is preserved while ``unnecessary'' geometric distortions are removed.
	
	On the target interval, a new map $F$ is defined to be piecewise linear with constant slope magnitude $\lambda>1$, where $\lambda = e^{h(f)}$. The lap structure of $f$ is reproduced by $F$, and the symbolic dynamics is matched exactly so that the relation $\phi \circ f = F \circ \phi$ is satisfied and $h(F) = h(f) = \log \lambda$ is obtained.
	
	\begin{definition}
		Let $f:I\to I$ be a continuous self-map of an interval $I $. A point $c\in I$ is called a \emph{critical point}\index{point!critical} of $f$ if $c$ is a local extremum point of $f$, that is, if $c$ admits an open neighbourhood in which $f(c)$ is a local maximum or minimum value. The \emph{critical set}\index{set!critical} of $f$ is the set of all local extrema:
		\[
		\operatorname{Crit}(f) := \{\, c \in I \mid c \ \text{is a critical point of} \ f \,\}.
		\]		
		The \emph{post-critical set}\index{set!post-critical} of $f$ is
		\[
		P(f)\;:=\;\overline{\bigcup_{n\geq 1} f^{n}\bigl(\operatorname{Crit}(f)\bigr)},
		\]
		where the overline denotes closure in $I$. 
		The map $f$ is called \emph{post-critically finite} (\emph{PCF})\index{map!post-critically finite} if $P(f)$ is a finite subset of $I$. 
	\end{definition}
	
	Theorem~\ref{thm:existence/of/uniform/expander} can be extended to the case of PCF maps. More precisely, if $f$ is a piecewise-monotone PCF map of an interval $I$, then there is a PCF uniform $\lambda$-expander $F$ semi-conjugate to $f$. Since entropy is invariant under topological conjugacy, we have $h(f) = h(F)$.
	
	\subsection{Incidence Matrices for Interval Maps.}
	\label{subsec:incidence-matrices-interval-maps}
	In this part, we review the construction of incidence matrices associated with one-dimensional dynamical systems.
	
	If \(f:I\to I\) is a continuous piecewise-monotone map, then a finite partition of \(I\) into intervals on which \(f\) is monotone allows one to associate a non-negative transition matrix to \(f\). In the simplest Markov situation, let	$I=I_1\cup\cdots\cup I_n$ be a finite interval partition such that \(f\) is monotone on each \(I_i\) and such that the image of each \(I_i\) is a union of elements of the partition. The corresponding incidence matrix \(A=(a_{i,j})\) is defined by
	\[
		a_{i,j}	=
		\begin{cases}
			1, & \text{if } f(I_i)\supseteq I_j,\\
			0, & \text{otherwise.}
		\end{cases}
	\]
	Thus \(a_{i,j}\) records only whether the interval \(I_j\) is covered by the image of \(I_i\). In this Markov case \(A\) is the adjacency matrix of the associated symbolic dynamics, and the entropy of the corresponding Markov system is
	\[
	h_{\mathrm{top}}(f)=\log\rho(A),
	\]
	where \(\rho(A)\) denotes the spectral radius of \(A\).
	
	For the more general class of piecewise-monotone maps with finitely many critical points, one uses the extended incidence matrix as given in the following:
	\begin{definition}
		Let \(V=\{v_0<\cdots<v_N\}\) be a finite forward-invariant set containing
		the critical values of \(f\), and set \(U_j=(v_{j-1},v_j)\). The
		\emph{extended incidence matrix}\index{matrix!extended incidence}
		\(\widehat A=(\widehat a_{ij})\) is defined by
		\[
		\widehat a_{ij}
		=
		\left|
		\left\{
		\text{connected components of } f^{-1}(U_j)\cap U_i
		\right\}
		\right|.
		\]
	\end{definition}
	Thus \(\widehat a_{ij}\) counts the number of monotone branches of \(f\) from \(U_i\) crossing \(U_j\). Unlike the \(0\)-\(1\) incidence matrix, the extended incidence matrix keeps track of multiplicities of coverings.
	
	The matrix \(\widehat A\) is non-negative. Moreover, because the intervals \(U_j\) are linearly ordered, the non-zero entries in each column form a consecutive block. If \(c(j)_t\) and \(c(j)_b\) denote, respectively, the row numbers of the first and last non-zero entries in the \(j\)-th column, then these blocks satisfy the connectedness condition
	\[
	|c(j)_t-c(j-1)_b|\leq 1 \quad\text{or}\quad	|c(j)_b-c(j-1)_t|\leq 1
	\]
	for each \(j=2,\ldots,N\). In other words, one can pass from a non-zero entry in the first column to a non-zero entry in the last column without crossing a zero entry. Finally, by \cite[Proposition~3.3]{thurston/entropy/in/dimension/one}, the entropy is recovered from the spectral radius of the extended incidence matrix:
	\[
		h_{\mathrm{top}}(f)=\log\rho(\widehat A).
	\]
	\subsection{Proof of Lind's theorem.}	
	\label{sec:thurston/lind/theorem}
	Let \(\lambda > 1\) be a Perron number. Thurston gave a proof of Lind's theorem via a convex–geometric and combinatorial encoding, by an explicit construction of a one–dimensional dynamical system with topological entropy \(\log \lambda\). To this end, we let 
	$V_\lambda \;=\; \QQ(\lambda) \otimes_{\QQ} \RR$
	be the real vector space obtained from \(\QQ(\lambda)\) by extension of scalars. This decomposes as the product of copies of \(\RR\) for the real embeddings (listed usually as $\sigma_{1}, \sigma_{2}, \ldots, \sigma_{r}$) and \(\CC\) for the complex embeddings (listed usually as $\sigma_{r+1}, \overline{\sigma_{r+1}}, \ldots, \sigma_{r+s}, \overline{\sigma_{r+s}},$) of the number field \(\QQ(\lambda)\). Indeed, the dimension of the $\RR$ vector space $V_{\lambda}$ is $\deg(\lambda)$; furthermore, if $r$ is the number of real embeddings and $2s$ is the number of complex embeddings of $\QQ(\lambda)$ into $\CC$, then $d = r + 2s$ and $V_{\lambda} \simeq \RR^{r} \times \CC^{s}$.  Multiplication by \(\lambda\), denoted by $m_{\lambda}$, acts linearly on \(V_\lambda\). If one chooses $\{1,\lambda,\ldots,\lambda^{d-1}\}$ as a basis, the matrix representing $m_{\lambda}$ becomes the companion matrix of $\lambda$.  In particular, $\lambda$ is the dominant eigenvalue of $m_{\lambda}$, see \ref{sec:dominant/eigenvalue}. However, this matrix may not be non-negative. So, we let \(K_\lambda \subset V_\lambda\) denote the \emph{dominant cone} of $m_{\lambda}$, that is, the cone in $V_{\lambda}$ consisting of those vectors whose projection to the \(\lambda\)–eigenspace $E_{\lambda}$ is larger in modulus than the projection to any other invariant subspace. Then, $m_{\lambda}(K_{\lambda}) \subset \operatorname{int}(K_\lambda)$. 
	
	The map $\iota:\QQ(\lambda)\longrightarrow V_\lambda$ sending every element \(\alpha\in\QQ(\lambda)\) to
	\[
		\sigma_1(\alpha)\otimes\cdots\otimes\sigma_d(\alpha) \in V_{\lambda}
	\]
	has dense image. We write $\mathcal O_\lambda:=\mathcal O_{\QQ(\lambda)}$ for the ring of integers of \(\QQ(\lambda)\), that is, the set of elements of \(\QQ(\lambda)\) which satisfy a monic polynomial in \(\ZZ[x]\). Furthermore, the image \(\iota(\mathcal O_\lambda)\) is a lattice in \(V_\lambda\). For details we refer to \cite[\S~1.12]{neukirch/ant} and to \cite{schoof/arakelov}. As a result of the density  of $\iota(\QQ(\lambda))$, one may find finitely many (the finiteness is a result of Gordon's lemma, see, for instance, David A. Cox, Lectures on toric varieties, Lecture 1. Proposition 1.11\footnote{https://dacox.people.amherst.edu/lectures/coxcimpa.pdf}) elements, say $\alpha_{1},\ldots,\alpha_{n}$, of $V_{\lambda} \cap \iota(\QQ(\lambda))$ so that the cone over these points, say $C_{\alpha_{1},\ldots,\alpha_{n}}$, satisfies $m_{\lambda}(K_{\lambda}) \subset C_{\alpha_{1},\ldots,\alpha_{n}}\subset \operatorname{int}(K_{\lambda})$. The intersection $\iota(\mathcal{O}_{\QQ(\lambda)}) \cap C_{\alpha_{1},\ldots,\alpha_{n}}$ is also finitely generated, say by $g_{1},\ldots,g_{m}$. Moreover, $\lambda g_{i}$ can also be expressed as a linear combination of the generators where each coefficient is non-negative. So, we obtain a non-negative $m \times m$ matrix, say $A$, having $\lambda$ as an eigenvalue. The projections of the vectors $\lambda g_{i}$ onto the eigenspace $E_{\lambda}$ are all ``positive'', producing a positive eigenvector for $A$, which proves Lind's theorem.
	
	\subsection{Realization for Powers of Perron Numbers}
	\label{subsec:realization-powers-perron}
	To ensure that the resulting piecewise-linear map is PCF, the underlying transition matrix must satisfy strict parity and matching conditions. Thurston adapts the classical framework of Lind's Theorem to construct an incidence matrix meeting these combinatorial constraints via the following multi-step algebraic machinery:
	
	\begin{itemize}
		\item \textbf{Semigroup Generators:} One selects a finite sequence of generators $g_{1}, \dots, g_{k}$ for a positive cone $K_{\lambda} \subset \mathcal{O}_{\lambda}$ such that every generator appears at least once. Their partial sums must surjectively map to the quotient classes of the finite abelian group $\mathcal{O}_{\lambda} / 2 \mathcal{O}_{\lambda}$ and yield a vanishing total sum $T = \sum g_i$, which guarantees that $T \equiv 0 \pmod{2\mathcal{O}_{\lambda}}$.
		
		\item \textbf{The Slim Cone Lemma:} Since multiplication by the Perron number $\lambda$ on the projective completion of the cone $K_{\lambda}$ acts as a contraction with a unique attracting fixed point, there exists a sufficiently large integer $N$ such that $\lambda^{N} K_{\lambda} \subseteq 3T + K_{\lambda}$. 
		
		\item \textbf{Construction of a Uniform $\lambda^{N}$-Expander:} A base interval of total length $T$ is subdivided into consecutive segments whose lengths correspond to the sequence of generators $g_{i_1}, \dots, g_{i_k}$ (as measured under the primary real embedding). By an inductive combinatorial matching argument, one constructs a Markov interval map where the image of each segment $g_{i}$ is uniquely determined by its class $\lambda^{N}g_{i} \pmod{2\mathcal{O}_{\lambda}}$. This guarantees that the resulting map acts as a uniform $\lambda^{N}$-expander governed by a valid, integral incidence matrix that preserves the PCF property.
	\end{itemize}
		
%	Given \(\lambda\) and the generating set \(G\), choose a finite sequence of generators \((g_{i_1},\dots,g_{i_k})\) such that:
%	\begin{itemize}
%		\item Every generator from \(G\) appears at least once.
%		\item The partial sums represent all residue classes modulo \(2\mathcal{O}_\lambda\).
%		\item The total sum \(T \equiv 0 \pmod{2\mathcal{O}_\lambda}\).
%	\end{itemize}
	\noindent By the contraction property in \(K_\lambda^{\mathbb{Q}}\), there exists \(N \gg 0\) such that
	\[
	\lambda^N \cdot S_\lambda \subset 3T + S_\lambda.
	\]
	Subdivide the unit interval into subintervals of lengths \(g_{i_1},\dots,g_{i_k}\) (measured in the real embedding where \(\lambda\) is real \(>1\)). Define the image of each subdivision point so that \(\lambda^N g_{i_j}\) lands at a subdivision point with the correct parity modulo \(2\mathcal{O}_\lambda\); the difference is an \emph{even}, strictly positive combination of generators, realized by zig-zag traversals across the corresponding target intervals. This produces a uniform \(\lambda^N\)-expander with an extended incidence matrix.
	The reduction from $\lambda^{N}$ to $\lambda$ can be obtained by two methods: 	
	\begin{itemize}
		\item \textbf{Circle method.} View the circle as subdivided into \(N\) arcs with lengths \(1, \lambda, \dots, \lambda^{N-1}\), cyclically permuted, with the \(\lambda^N\)-expander inserted on the last arc and rescaled appropriately. The resulting map has constant slope \(\lambda\).
		\item \textbf{Interval method.} Take a quadratic interval map from the period–doubling cascade with entropy \(0\) and critical point of period \(N=2^n\). \emph{Blow up} the forward and backward orbits of the critical point into small intervals with summable lengths, and implant the \(\lambda^N\)-expander as the first–return map on the blown–up critical interval. The resulting piecewise-linear map is a uniform \(\lambda\)-expander.
	\end{itemize}
	
	\begin{remark}
		This construction can immediately be used to construct expanders for \emph{weak Perron} integers, which are defined as those real algebraic integers \(\lambda\)  for which \(\lambda \ge |\sigma(\lambda)|\) for all embeddings \(\sigma\). Indeed, in such a case all conjugates of $\lambda$ is an element of the closed disk of radius $|\lambda|$. Those that are on the boundary must be of the form $\zeta \lambda$ for some root of unity $\zeta$, i.e. $|\zeta| = 1$. Choose an integer $k$ so that for every such conjugate $\zeta^{k} = 1$. Then, $\lambda^{k}$ is Perron. The preceding construction leads to a \(\lambda^k\)-expander. Then the implant method can be used to yield a uniform \(\lambda\)-expander.
	\end{remark}

	\section{Entropies of surface homeomorphisms}
	\label{sec:entropies/of/surface/homeomorphisms}
	So far, we have seen that the dynamics of one-dimensional interval maps is fundamentally combinatorial and has a significant algebraic component. The passage to surfaces introduces new topological constraints. The goal of this section is to explain how the entropy of a surface homeomorphism, especially in the pseudo-Anosov case, is encoded by invariant measured foliations, train tracks, and the spectral radii of non-negative transition matrices.
	
	We first recall the Nielsen--Thurston classification. We then explain how train tracks convert pseudo-Anosov dynamics into Perron--Frobenius data. Finally, we discuss the resulting algebraic restrictions on dilatations, emphasizing the contrast with the interval case: while post-critically finite interval maps realize weak Perron numbers, pseudo-Anosov surface homeomorphisms satisfy stronger restrictions coming from invertibility, homology, and the intersection form. We refer to \cite{farb/margalit/MCG} for background.
	
	\subsection{The Nielsen-Thurston Classification}
	
	Let $S = S_{g,n}$ be a compact connected orientable surface of genus $g$ with $n$ punctures (or marked points) such that its Euler characteristic $\chi(S) = 2 -2g - n < 0$. Let $\mathrm{Mod}(S) := \pi_0(\mathrm{Homeo}^+(S))$ denote the mapping class group\index{mapping class group} of $S$. 
	Before stating the classification theorem, we recall the main concepts that appear in what follows. A \emph{simple closed curve}\index{curve!simple closed} on \(S\) is an embedded copy of the circle, that is a continuous embedding
	\[
		\gamma:S^1\hookrightarrow S.
	\]
	As always, we identify the function with the image $\gamma(S)$. Notice that, the curve has no self-intersections. A finite collection
	\[
		\mathcal C=\{\gamma_1,\ldots,\gamma_k\}
	\]
	of simple closed curves is called a \emph{multicurve}\index{curve!multi} if its components are pairwise disjoint and satisfy certain non-degeneracy conditions. More precisely, in the Nielsen--Thurston classification one considers finite systems of disjoint, essential, non-parallel simple closed curves. Let us explain each of
	these terms.
	
	\begin{definition}[Finite system]
		A \emph{finite system}\index{curve system!finite} of simple closed curves on \(S\) is a finite set
		\[
		\mathcal C=\{\gamma_1,\ldots,\gamma_k\}
		\]
		where each \(\gamma_i\) is a simple closed curve on \(S\). The word \emph{finite} means that only finitely many isotopy classes of curves are involved.
	\end{definition}
	
	\begin{definition}[Disjoint]
		The system \(\mathcal C\) is called \emph{disjoint}\index{curve system!disjoint} if the curves can be chosen so that
		\[
		\gamma_i\cap \gamma_j=\varnothing
		\qquad\text{for all } i\neq j.
		\]
		Equivalently, the curves are represented by pairwise disjoint embedded circles on the surface.
	\end{definition}
	
	\begin{definition}[Essential]
		A simple closed curve \(\gamma\subset S\) is called \emph{essential}\index{curve!essential} if it is not null-homotopic and not homotopic into a boundary component or into a puncture. Equivalently, \(\gamma\) is essential if it does not bound a disk, a once-punctured disk, or an annulus whose other boundary component is a boundary component of \(S\).
	\end{definition}
	
	\begin{definition}[Non-parallel]
		Two essential simple closed curves \(\gamma_1,\gamma_2\subset S\) are called \emph{parallel} if they are isotopic, or equivalently if they cobound an embedded annulus in \(S\). A system \(\mathcal C\) is called \emph{non-parallel} if no two distinct components of \(\mathcal C\) are parallel.
	\end{definition}
	Thus a finite system of disjoint, essential, non-parallel simple closed curves is a finite multicurve
	\[
	\mathcal C=\{\gamma_1,\ldots,\gamma_k\}
	\]
no two of whose components represent the same isotopy class. Such multicurves are precisely the curves along which a reducible mapping class can be decomposed.

	The  Nielsen--Thurston classification theorem, see for instance \cite{Thurston1988GeometryDynamics}, or \cite{FathiLaudenbachPoenaru2012ThurstonWork,CassonBleiler1988Automorphisms}, asserts that every mapping class $[f] \in \mathrm{Mod}(S)$ contains a canonical representative $\phi \in [f]$ satisfying exactly one of the following geometric characterizations:
	\begin{enumerate}
		\item \textbf{Periodic:}\index{mapping class!periodic}There exists an integer $k \ge 1$ such that $\phi^k = \mathrm{id}_S$. In this setting, the topological entropy satisfies $h(\phi) = 0$. In this case, by Kerckhoff's resolution of the Nielsen realization problem, \cite{Kerckhoff1983NielsenRealization}, a periodic map can be realized as a strict isometry (a rotation or reflection) under some hyperbolic, flat, or spherical metric on the surface.
		
		\item \textbf{Reducible:}\index{mapping class!reducible} There exists a non-empty, finite system of disjoint, essential, non-parallel simple closed curves $\Gamma = \{\gamma_1, \ldots, \gamma_m\}$ on $S$ that is invariant under $\phi$ (i.e., $\phi(\Gamma) = \Gamma$). The dynamics can then be analyzed inductively by cutting $S$ along $\Gamma$ into subsurfaces of strictly smaller complexity. Geometrically, if one cuts the surface along these invariant curves, the map $\phi$ permutes or maps the resulting smaller component pieces into themselves. On those smaller pieces, the map will again be either periodic or pseudo-Anosov.
		
		\item \textbf{Pseudo-Anosov:}\index{mapping class!pseudo-Anosov} There exists a pair of transverse, singular measured foliations on $S$, namely the stable foliation $(\mathcal{F}^s, \mu_s)$ and the unstable foliation $(\mathcal{F}^u, \mu_u)$, and a real constant $\lambda > 1$, known as the \emph{dilatation} or \emph{stretch factor}, such that 
		\[
		\phi(\mathcal{F}^u, \mu_u) = (\mathcal{F}^u, \lambda \mu_u) \quad \text{and} \quad \phi(\mathcal{F}^s, \mu_s) = (\mathcal{F}^s, \lambda^{-1} \mu_s).
		\]
		That is, under the map $\phi$, paths along the stable foliation are contracted by a factor of $1/ \lambda$, while paths along the unstable foliation are stretched by a factor of $\lambda$.
	\end{enumerate}
	
	In the pseudo-Anosov world, the continuous stretching of the unstable foliation by $\lambda$ and the simultaneous contraction of the stable foliation by $\lambda^{-1}$ govern the global mixing properties of the map. This explicit geometric scaling directly determines the asymptotic orbit complexity, yielding an exact formulation for the topological entropy:
	\[
	h(\phi) = \log(\lambda).
	\]
	Furthermore, a theorem of Thurston states that $\log(\lambda)$ minimizes the topological entropy within the entire mapping class $[f]$, making the pseudo-Anosov dilatation a fundamental topological invariant of the class.
	
	\subsection{Train Tracks and Incidence Matrices}
	
	To bridge the gap between the continuous geometry of measured foliations and the discrete spectrum of linear operators (and hence combinatorics), we use the theory of \emph{train tracks}. This viewpoint allows us to translate the continuous topological dynamics of a pseudo-Anosov homeomorphism into concrete linear algebra. In particular, we obtain a finite-dimensional projection of the action on measured foliations and therefore a direct link to Perron--Frobenius theory.
	
	\begin{definition}[\cite{Thurston2023geometry/of/3/manifolds,penner-train-tracks}]
		A \emph{train track}\index{train track} \(\tau\) on a compact smooth surface \(S\) is an embedded multigraph, that is, a finite graph embedded in \(S\) in which multiple edges between the same pair of vertices and loops are allowed, satisfying the following conditions:
		\begin{enumerate}
			\item Each vertex $v$ (called a \emph{switch}) is assigned a unique tangent line in $T_v S$. The edges (called \emph{branches}) incident to $v$ are tangent to this line. The edges are separated into two sets: those entering from the left, $E_{in}(v)$, and those leaving from the right, $E_{out}(v)$.
			\item For every connected component $U$ of $S \setminus \tau$, the completion of $U$ with respect to a path metric has a negative generalized Euler characteristic, defined by
			\[
			\chi'(U) = \chi(U) - \frac{1}{2}\sum_{i} c_i < 0,
			\]
			where $c_i$ are the inward-pointing cusps on the boundary of $U$. 
		\end{enumerate}
	\end{definition}
	
	Remark that the last condition eliminates problematic components such as monogons without punctures, digons, or annuli. A train track $\tau$ acts as a template for families of curves and measured foliations\index{measured foliation}. A regular neighborhood $N(\tau)$ can be viewed as a union of thin bands (associated with branches) meeting at junctions (associated with switches). A measured foliation $\mathcal{F}$ is said to be \emph{carried}\index{measured foliation!carried} by $\tau$ (denoted $\mathcal{F} \prec \tau$) if it can be embedded inside $N(\tau)$ such that its leaves run parallel to the branches. 
	
	To formalize this algebraically, we assign a non-negative real weight $w_i = w(e_i)$ to each branch $e_i \in E(\tau) = \{e_1, \ldots, e_k\}$. For the weights to yield a well-defined transverse measure under the carrying relation, they must satisfy conservation of flow at every switch. This leads to the \emph{switch conditions}\index{switch condition} (analogous to Kirchhoff's laws that appear in circuit analysis):
	\begin{equation}\label{eq:switch_condition}
		\sum_{e_i \in E_{in}(v)} w_i = \sum_{e_j \in E_{out}(v)} w_j \quad \forall v \in V(\tau).
	\end{equation}
	The set of all vectors $w \in \RR^k_{\ge 0}$ satisfying \eqref{eq:switch_condition} forms a convex cone, denoted $\mathcal{W}(\tau)$, which parametrizes the space of measured foliations carried by $\tau$.
	
	Now let $\phi: S \to S$ be a pseudo-Anosov homeomorphism. By choosing a suitable train track $\tau \subset S$, we can arrange for $\tau$ to be \emph{invariant} under $\phi$, meaning that the image graph $\phi(\tau)$ can be isotoped into the regular neighborhood $N(\tau)$ without back-tracking. This geometric inclusion allows us to define the combinatorial \emph{incidence (or transition) matrix} $M = (m_{i,j}) \in \mathrm{Mat}_{k \times k}(\ZZ_{\ge 0})$, where
	\[
	m_{i,j} = \text{the number of times } \phi(e_j) \text{ traverses the branch } e_i.
	\]
	In addition, as $\phi$ is an orientation-preserving homeomorphism that maps the smooth structure of the switches properly, the linear transformation defined by $M$ preserves the switch relations. Thus, $M(\mathcal{W}(\tau)) \subset \mathcal{W}(\tau)$.
	
	The pseudo-Anosov map's global mixing properties imply specific algebraic properties of the matrix $M$. Indeed, as a pseudo-Anosov map minimizes entropy and possesses dense leaves in its invariant foliations, the track cannot contain any invariant sub-tracks or periodic loops that do not interact. Topologically, every branch eventually expands under iteration to slide across every other branch. Algebraically, this implies that the incidence matrix $M$ is \emph{primitive} (compare with Sec~\ref{sec:dominant/eigenvalue}): there exists an integer $\nu \ge 1$ such that every entry of $M^\nu$ is strictly positive.
	
	By the Perron--Frobenius theorem (Theorem~\ref{thm:frobenius/theorem/for/irreducible/matrices}), the primitivity of $M$ guarantees that:
	\begin{enumerate}
		\item $M$ has a unique, simple, real dominant eigenvalue $\rho(M) > 1$ (the Perron root) equal to its spectral radius.
		\item The corresponding unique right eigenvector $v_{\mathrm{PF}}$ and left eigenvector $u_{\mathrm{PF}}$ can be chosen to have strictly positive entries.
	\end{enumerate}
	Since $M$ maps the subspace defined by the linear system \eqref{eq:switch_condition} into itself, and $v_{\mathrm{PF}}$ is unique up to scaling, $v_{\mathrm{PF}}$ must lie within the positive cone $\mathcal{W}(\tau)$. The coordinates of $v_{\mathrm{PF}}$ define the precise transverse measures on the branches required to reconstruct the unstable foliation $\mathcal{F}^u$. Under the action of $\phi$, these weights scale uniformly by $\rho(M)$. Thus, the combinatorial spectral radius exactly recovers the geometric dilatation:
	\[
	\rho(M) = \lambda, \quad \text{and consequently,} \quad h(\phi) = \log(\rho(M)).
	\]
	
	\subsubsection{An Explicit Example: surface of genus $2$}\label{example}
	
\begin{figure}
\includegraphics[width=\textwidth]{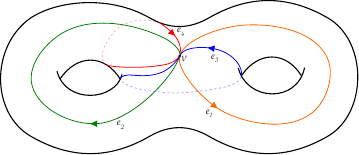}
\caption{A genus 2 surface with a train track invariant under a pseudo-Anosov map.}
\label{fig:genus2}
\end{figure}

	On a surface of genus $2$, $S$, consider the train track $\tau \subset S$ containing a single switch $v$ and four edges or branches $E(\tau) = \{e_1, e_2, e_3, e_4\}$. The edges $e_1, e_2, e_3$ and $e_4$ emanate from $v$, and smoothly merge to terminate at $v$ in the order $e_3, e_1, e_4, e_2$, see Figure~\ref{fig:genus2}. The complement of $\tau$ in $S$ is a single, simply-connected octagon. There is no non-trivial weight relation, $\mathcal{W}(\tau)\cong\mathbb{R}^4$.
	
	One can define a pseudo-Anosov mapping class $\phi: S \to S$ whose topological action on $\tau$ sends $e_i$ along the concatenated paths:
	\begin{itemize}
	\item $e_1 \mapsto e_1 \cdot e_2 \cdot e_3 \cdot e_1 \cdot e_3$,
	\item $e_2 \mapsto e_1 \cdot e_2 \cdot e_3 \cdot e_2 \cdot e_3 \cdot e_1 \cdot e_2  \cdot e_4 \cdot e_3$,
	\item $e_3 \mapsto e_1 \cdot e_2 \cdot e_3 \cdot e_2 \cdot e_3 \cdot e_1 \cdot e_3$,
	\item $e_4 \mapsto e_1 \cdot e_2 \cdot e_4 \cdot e_3$.	
\end{itemize}

 Counting the occurrences of each branch yields the following non-negative incidence matrix $M$:
	\[
	M = \begin{pmatrix} 
	 2 & 2 & 2 & 1 \\
	 1 & 3 & 2 & 1 \\
	 2 & 3 & 3 & 1 \\
	 0 & 1 & 0 & 1 \\ 
	 \end{pmatrix}.
	\]
	
	Note that
	\[
	M^2 =  \begin{pmatrix}  10 & 17 & 14 & 7 \\
 9 & 18 & 14 & 7 \\
 13 & 23 & 19 & 9 \\
 1 & 4 & 2 & 2  \end{pmatrix},
	\]

	that is, $M^{2}$ is positive, meaning $M$ is primitive, and the Perron--Frobenius theorem applies directly. The eigenvalues are roots of the characteristic polynomial
	\[
	\det(M - xI) = x^4-9 x^3+16 x^2-9 x+1 = (x-1)^2 \left(x^2-7 x+1\right)
	\]
	The roots are $x=1$ and $x = (7 \pm 3\sqrt{5})/2$. The spectral radius is dictated by the unique dominant positive root: $\lambda = \rho(M) = (7 + 3\sqrt{5})/2$. To find the invariant measure of the unstable foliation, we compute the right Perron--Frobenius eigenvector $v_{\mathrm{PF}}$ associated with $\lambda$:

	\[
	v_{\mathrm{PF}} = \begin{pmatrix} 7+3 \sqrt{5} \\ 2 \\ 2 \\ 7-3 \sqrt{5} \end{pmatrix}.
	\]
	Note that $(7-3\sqrt{5})/2 = 1/\lambda > 0$, so all entries are strictly positive as dictated by Perron--Frobenius theory. 

	\subsection{Algebraic Restrictions on the Dilatation}
	
	There are fundamental differences between the arithmetic nature of dilatations in one-dimensional dynamics and in surface dynamics. Recall from Section~\ref{sec:thurston/lind/theorem} that, for post-critically finite interval maps, \(\exp(h_{\mathrm{top}})\) can be any weak Perron number. For pseudo-Anosov surface homeomorphisms, the train-track construction still gives a non-negative integral transition matrix, and therefore the dilatation \(\lambda\) is a Perron--Frobenius eigenvalue. In particular, \(\lambda\) is a Perron number. However, the fact that the underlying map is an invertible surface homeomorphism imposes additional arithmetic restrictions.
	
	First, pseudo-Anosov dilatations are algebraic units. Equivalently, if 
	\(K=\QQ(\lambda)\), then
	\[
		\lambda\in \mathcal O_K^\times .
	\]
	This is a stronger condition than being the spectral radius of an arbitrary non-negative integral matrix: such a matrix need not be invertible over \(\ZZ\), whereas a pseudo-Anosov map comes from an invertible topological dynamical system.
	
	There is a further restriction in the orientable case. Suppose that the stable and unstable invariant foliations of a pseudo-Anosov homeomorphism \(\phi:S\to S\) are orientable. Then their transverse measures determine cohomology classes which are eigenvectors for 
	\[
		\phi^*:H^1(S;\RR)\longrightarrow H^1(S;\RR)
	\]
	with eigenvalues \(\lambda\) and \(\lambda^{-1}\). Since \(\phi\) is orientation-preserving, \(\phi^*\) preserves the algebraic intersection form on \(H^1(S;\ZZ)\). Hence the characteristic polynomial of \(\phi^*\) is reciprocal, see \cite{rafiqi2024symplectic} for a further study. Consequently, in this orientable case, the factor containing \(\lambda\) is reciprocal, and \(\lambda^{-1}\) is a Galois conjugate of \(\lambda\). Without the orientability hypothesis, the dilatation need not occur as an eigenvalue of the homological action. Thus one should not conclude, in general, that the minimal polynomial of \(\lambda\) is reciprocal. The correct general restriction for closed surfaces is Fried's theorem \cite{fried1985growth}: the dilatation of a pseudo-Anosov homeomorphism on a closed surface is a bi-Perron unit.
	
	\begin{definition}
		A real algebraic unit \(\lambda>1\) is called a	\emph{bi-Perron unit}\index{algebraic number!bi-Perron} if every Galois conjugate \(\mu\) of \(\lambda\), (except possibly one of $\pm1/\lambda$), satisfies
		\[
			\lambda^{-1}< |\mu|< \lambda .
		\]
	\end{definition}
	Thus pseudo-Anosov dilatations on closed surfaces form a much smaller class than the weak Perron numbers appearing in interval dynamics. The converse problem remains subtle: determining which bi-Perron units, or  which non-negative transition matrices, are realized by pseudo-Anosov homeomorphisms is still a deep problem. Explicit sufficient conditions on aperiodic incidence matrices realizing pseudo-Anosov dilatations were given in \cite{baik2016constructing}.

	In summary, the interplay between topological constraints and dynamical complexity yields a rigid hierarchy of algebraic classifications:
	\[
	{\small
	\left\{ \begin{matrix} \text{pseudo-Anosov} \\ \text{dilatations on} \\ \text{closed surfaces} \end{matrix} \right\} \subseteq_? \left\{ \begin{matrix}
		\text{bi-Perron} \\ \text{units}
	\end{matrix} \right\} \subsetneq \left\{ \begin{matrix}
		\text{Perron} \\ \text{numbers}
	\end{matrix}\right\} \subsetneq \left\{\begin{matrix}
	\text{weak Perron} \\ \text{numbers}
	\end{matrix}\right\}}
	\]
	The strictness of these nested algebraic sets reflects the increasing density of geometric data imposed on the underlying system as we move from general one-dimensional maps to invertible, closed surface dynamics.
	
	\subsubsection{Stretch-factors on infinite type surfaces} 
	Interestingly though, even in the case of surface homeomorphisms, removing the compactness condition does realize all weak-Perron numbers as stretch-factors, as was recently shown in \cite{chenxi2026}:
	
	\begin{theorem}[Hillen-Loving-Wu]
	A number $\lambda$ is the stretch factor of an end-periodic homeomorphism of a connected infinite-type surface with finitely many ends if and only if $\lambda$ is a weak Perron number.
	\end{theorem}
	
	\emph{End-periodic} homeomorphisms are defined on surfaces of infinite type, and have many similarities with Thurston's pseudo-Anosov maps of compact surfaces. In the examples constructed by \cite{chenxi2026}, the non-compact, non-planar, infinite-type surfaces do have transverse stable and unstable measured singular foliations, and although the number of singularities may be infinite, there are only finitely many accumulation points of the singularities. One foliation is stretched, and the other contracted, by an algebraic integer $\lambda\geq1$, called the \emph{Handel-Miller stretch factor}. 
		
	\subsection{Generalization to Graph Maps and Free Group Automorphisms}
	
	Thurston's paradigm of extracting discrete matrix dynamics from continuous topological maps extends naturally beyond surfaces to the purely algebraic setting of free groups. Let $F_n$ be a free group of rank $n \ge 2$, and let $\mathrm{Out}(F_n) = \mathrm{Aut}(F_n) / \mathrm{Inn}(F_n)$ denote its outer automorphism\index{outer automorphism} group. The study of $\mathrm{Out}(F_n)$ shares a deep, structural dictionary with the mapping class group $\mathrm{Mod}(S)$. 
	
	Just as a mapping class is realized as a homeomorphism of a surface, an element $\Phi \in \mathrm{Out}(F_n)$ can be topologically realized as a self-homotopy-equivalence $f \colon G \to G$ of a finite connected metric graph $G$ with fundamental group $\pi_1(G) \simeq F_n$. To translate the dynamical complexity of $\Phi$ into linear algebra, Bestvina and Handel, \cite{bestvina1992train}, developed the theory of \emph{graph train tracks}, providing an exact combinatorial analogue to surface train tracks.
	
	\begin{definition}
		Let $G$ be a finite graph with no valence-one vertices. A continuous map $f \colon G \to G$ is a \emph{topological representative}\index{outer automorphism!topological representative of} of $\Phi \in \mathrm{Out}(F_n)$ if $f$ sends vertices to vertices, edges to paths, and induces $\Phi$ on the fundamental group. 
		Such a representative is called a \emph{graph train track map}\index{train track!graph} if for every edge $e \in E(G)$ and every integer $m \ge 1$, the path $f^m(e)$ is locally injective (i.e., it contains no non-trivial backtracking or cancellations).
	\end{definition}
	
%	For a graph train track map $f$, the action on the edges defines a non-negative integral incidence matrix $M \in \mathrm{Mat}_{k \times k}(\ZZ_{\ge 0})$, where $k = |E(G)|$. The entry $m_{i,j}$ counts the number of times the path $f(e_j)$ traverses the edge $e_i$ in either direction. Because there is no cancellation upon iteration, the word length of paths grows exactly according to the powers of $M$. Thus, if the action is ergodic (meaning $M$ is primitive), the topological entropy of the map is given precisely by $h(f) = \log(\lambda)$, where $\lambda = \rho(M)$ is the Perron--Frobenius eigenvalue.
	
	For a graph train track map \(f:G\to G\), let \(M_f\) be its transition matrix, whose \((i,j)\)-entry counts the number of occurrences of the edge \(e_i\) in the tightened edge-path \(f(e_j)\), ignoring orientation. The \emph{expansion constant} of \(f\) is defined to be
	$\lambda(f):=\rho(M_f),$ where \(\rho(M_f)\) denotes the spectral radius of \(M_f\).\footnote{In the pseudo-Anosov setting, the number \(\lambda>1\) is more commonly called the \emph{stretch factor} or \emph{dilatation}. In the graph train-track setting, especially in formulations of Thurston's entropy theorem, the same type of quantity is often called the \emph{expansion constant}. In both cases it is the exponential of the topological entropy and, when a transition matrix \(M\) is present, it is given by \(\lambda=\rho(M)\). Thus the terms \emph{stretch factor} and \emph{expansion constant} refer here to the same Perron--Frobenius growth rate.} Equivalently, \(\lambda(f)\) is the exponential growth rate of edge lengths:
	\[
	\lambda(f) = \lim_{n\to\infty} \left|f^n(e)\right|^{1/n}
	\]
	for any edge \(e\) whose iterates lie in the exponentially growing stratum. In the irreducible/ergodic case this limit is independent of \(e\). The topological entropy is therefore
	\[
	h_{\mathrm{top}}(f)=\log \lambda(f).
	\]
	Thus, if \(h=h_{\mathrm{top}}(f)\), the corresponding expansion constant is \(\lambda=\exp(h)\).
	
	In this context, Thurston established a beautiful and definitive analogue to his interval map realization theorem, bridging the topology of free group automorphisms directly back to the arithmetic of weak Perron numbers:
	
	\begin{theorem}[Thurston \cite{thurston/entropy/in/dimension/one}] \label{thm:thurston_outFn}
		A positive real number \(h\) occurs as the topological entropy of an ergodic graph train track representative of an outer automorphism \(\Phi\in \operatorname{Out}(F_n)\) if and only if $\lambda:=\exp(h) $ is a weak Perron number.
	\end{theorem}
	Here, we use the phrase \emph{ergodic graph train track representative} to mean a graph train track representative \(f:G\to G\) whose transition matrix \(M_f\) is primitive. Equivalently, there exists \(N\geq 1\) such that every edge of \(G\) occurs in the tightened path \(f^N(e)\) for every edge \(e\) of \(G\). This condition rules out invariant proper subgraphs and periodic decompositions of the edge dynamics. Under this hypothesis, Perron--Frobenius theory gives a unique dominant eigenvalue \(\lambda(f)=\rho(M_f)>1\), and this number is the expansion constant of \(f\).
		
	\subsubsection{Topological Obstructions versus Algebraic Flexibility}
	
	Theorem~\ref{thm:thurston_outFn} completes the loop of Thurston's contributions to Perron--Frobenius theory, but it introduces a vital concluding caveat that distinguishes one-dimensional algebraic topology from two-dimensional geometry. 
	
	Recall that for a pseudo-Anosov homeomorphism of a closed surface, the global invertibility of the surface map imposes restrictions stronger than those for arbitrary non-negative transition matrices. In particular, Fried's theorem implies that its dilatation is a bi-Perron unit. Under additional orientability hypotheses, the dilatation is detected by the homological action, whose characteristic polynomial is reciprocal.
	
	In contrast, a graph map $f \colon G \to G$ representing an outer automorphism is only a ``homotopy equivalence'', not a homeomorphism. This means that while the induced map on homology $f_* \colon H_1(G; \ZZ) \to H_1(G; \ZZ)$ must be an invertible element of $\mathrm{GL}_n(\ZZ)$ (and thus its eigenvalues are algebraic units), the action on the raw edge space $\RR^E$ is far less constrained. 
	
	Cellular operations on graphs allow edges to stretch, fold onto each other, or collapse entirely in a locally non-invertible fashion. The dimension of the incidence matrix $M$ is the number of edges $k$, which is determined by the Euler characteristic $\chi(G) = 1-n$. Thus, $k = |V| + n - 1$. Whenever $|V| > 1$, the dimension of the transition matrix strictly exceeds the homological dimension $n$. The ``extra'' eigenvalues of $M$ are associated with the combinatorics of vertices and folds rather than the invertible topology of the group.
		
	Consequently, the incidence matrix of a graph train track does not inherit any intrinsic reciprocal or unit constraints. The expansion constant $\exp(h)$ for a general free group automorphism does not have to be an algebraic unit, therefore one expects full algebraic flexibility in this case.
	
	\subsubsection{The Perron--Frobenius Degree of an Expansion Constant}
	
	Theorem~\ref{thm:thurston_outFn} guarantees that every weak Perron number can be realized as the spectral radius of some primitive integral matrix $M$, which corresponds to the transition matrix of a graph train track. A natural measure of the topological complexity of this realization is the minimal dimension of such a matrix.
	\begin{definition}
		Let \(\lambda\) be a Perron number. The \emph{Perron--Frobenius degree}\index{degree!Perron--Frobenius} of \(\lambda\), denoted by \(\deg_{\mathrm{PF}}(\lambda)\), is the smallest integer \(n\) for which there exists an \(n\times n\) non-negative integral primitive matrix whose spectral radius is \(\lambda\).
	\end{definition}
	 In a similar fashion, the \emph{irreducible Perron--Frobenius degree}\index{degree!irreducible Perron--Frobenius}, denoted by $\deg_{\mathrm{irrPF}}(\lambda)$, is defined as the smallest possible dimension of a non-negative, integral, irreducible matrix with spectral radius $\lambda$. 
	
	Because cellular graph operations introduce non-invertible, local folding, these matrix dimensions are generally much larger than the algebraic degree $d = [\QQ(\lambda): \QQ]$ of the expansion constant. The quantitative relationship between these degrees is highly non-trivial and computing them directly is difficult.
	
	By analyzing the layout of the Galois conjugates in the complex plane, Yazdi~\cite[Theorem~1.1]{yazdi2020lower} provided a geometric lower bound for the {P}erron--{F}robenius degree of a {P}erron number $p$ that has at least one non-real Galois conjugate $p'$, such that \[\tan^{-1}\left(\frac{p-\text{Re}(p')}{|\text{Im}(p')|}\right)\leq1.\] The left side of the above inequality, call it $\eta$, represents an angle determined by the geometric layout of these conjugates in the complex plane. If $\eta\leq1$, the Perron-Frobenius degree of $p$ is bounded below by:
	\[
	\deg_{\mathrm{PF}}(\lambda) \geq \frac{2\pi}{3 \eta},
	\]
	A striking application of this result is the proof that there exist cubic Perron numbers (and bi-Perron numbers of degree at most $6$) whose Perron–Frobenius degrees are arbitrarily large. Dynamically, this demonstrates that a mapping class with a very simple algebraic dilatation may necessitate an arbitrarily complex combinatorial train track to be topologically realized.
	
	The analogous question of explicit upper bounds has also been successfully answered, expressing the maximum necessary dimension in terms of deep number-theoretical data. By modifying Thurston’s proof of Lind’s realization theorem and introducing a methodology to control the number of integral points inside the invariant cone, Yazdi~\cite[Theorem~1.3]{yazdi2022nonnegative} established explicit upper bounds for $\deg_{\mathrm{irrPF}}(\lambda)$. These bounds depend on the algebraic degree $d$, the discriminant of the number field $D_{K}$ for $K = \QQ(\lambda)$, the minimal Hermite-like thickness of the ring of integers $\tau_{\min}(\mathcal{O}_{K})$, and the spectral ratio $\rho$.\footnote{The spectral ratio\index{algebraic number!spectral ratio of} $\rho$ of an algebraic number is defined as the modulus of its largest Galois conjugate strictly less than $\lambda$, divided by $\lambda$.} The upper bounds are given by:
	\begin{itemize}
		\item $\deg_{\mathrm{irrPF}}(\lambda) \leq \left(\frac{8d}{1 - \rho}\right)^{d^{2}}\tau_{\min}(\mathcal{O}_{K})^{d/2}$, and 
		\item $\deg_{\mathrm{irrPF}}(\lambda) \leq \left( \frac{8d}{1-\rho}\right)^{d^{2}} \sqrt{D_{K}}$.
	\end{itemize}
	The minimum of these two values yields a comprehensive upper bound for the irreducible case, which we denote by $B(\QQ(\lambda),\rho)$.
	
	These quantitative upper bounds can be transitioned from irreducible matrices to strictly primitive matrices by utilizing a structural relationship between the two types and applying a ``shift'' in the spectral radius. Suppose $A$ is an integral irreducible matrix with spectral radius $\mu$. By setting $\lambda = 1 + \mu$, one observes that the matrix $I + A$ will have spectral radius $\lambda$. In this context, two key observations apply:
	\begin{itemize}
		\item The matrix $I + A$ is naturally irreducible and has strictly positive diagonal entries, rendering it primitive. 
		\item If $\lambda > 1 + \frac{4}{1-\rho}$, then $\lambda - 1$ is also a Perron number generating the same number field, $\QQ(\lambda) = \QQ(\lambda-1)$. Thus, the upper bounds calculated for the irreducible matrix $A$ remain valid. 
	\end{itemize}
	
	By bounding the finite number of exceptions where $\lambda$ is too small to apply this matrix shift, this methodology yields a final, explicit upper bound for the primitive degree. Indeed, because the degree $d$ of the number field $\QQ(\lambda)$ is fixed and the moduli of all Galois conjugates are constrained by the spectral ratio $\rho$, classical finiteness results in algebraic number theory (such as Northcott's theorem, see for instance \cite{hindry2000diophantine}, on sets of algebraic integers with bounded height) guarantee that the subset of Perron numbers in $\mathbb{Q}(\lambda)$ intersecting the bounded interval $[1, 1 + \frac{4}{1-\rho}]$ is finite.
	
	One can define $\kappa(\QQ(\lambda),\rho)$ as the maximum of the Perron--Frobenius degrees among all these exceptional Perron numbers. By taking the maximum of this exceptional bound and the shifted bound from the irreducible case, one arrives at the definitive upper bound for the primitive degree:
	\[ 
	\deg_{\mathrm{PF}}(\lambda) \leq \max\Big\{2d^{2} B(\QQ(\lambda),\rho), \ \kappa(\QQ(\lambda),\rho)\Big\},
	\]
	where $B(\QQ(\lambda),\rho)$ represents the previously determined upper bound for the irreducible case. Together, these bounds frame the exact maximum size of the incidence matrix—and thus the combinatorial complexity of the graph train track—needed to topologically realize a given algebraic expansion constant. The transition from the irreducible to the strictly primitive case, alongside these explicit bounds, was established by M.~Yazdi \cite[Theorem~1.3]{yazdi2022nonnegative}.
	
	\section{Least stretch factors}
	
	A fundamental problem in the quantitative study of mapping class groups is understanding the minimal dynamic complexity of a pseudo-Anosov homeomorphism on a given surface. Because the set of stretch factors on a surface of fixed topological type is discrete, there exists a unique minimal dilatation $\delta_g > 1$ among all pseudo-Anosov mapping classes on the closed orientable surface $S_g$. 
	
	The minimal dilatation $\delta_g$ carries information about the geometry of the moduli space of Riemann surfaces, $\mathcal{M}_g$. In the Teichm\"uller metric, a pseudo-Anosov conjugacy class corresponds to a closed geodesic, and its length is precisely $\log(\lambda)$. Thus, determining $\delta_g$ is equivalent to finding the length of the shortest closed geodesic (the systole) in $\mathcal{M}_g$. Consequently, the study of  entropies for surface homeomorphisms—realized as Perron numbers—occupies a central role in this context.
	
	\subsection{Asymptotic Bounds and Exact Values}
	
	A fundamental result concerning the asymptotic behavior of \(\delta_g\) as
	\(g\to\infty\) is due to Penner~\cite{penner1991bounds}. More precisely,
	for the closed surface \(S=S_g\), Penner established the explicit bounds
	\[
	\frac{\log 2}{12g-12}
	\leq
	\log\delta_g
	\leq
	\frac{\log 11}{g}.
	\]
	Since \(12g-12\leq 12g\), we obtain
	\[
	\frac{\log 2}{12g}
	\leq
	\frac{\log 2}{12g-12}
	\leq
	\log\delta_g
	\leq
	\frac{\log 11}{g}.
	\]
	Thus
	\[
	\frac{\log 2}{12}\cdot\frac{1}{g}
	\leq
	\log\delta_g
	\leq
	\log 11\cdot\frac{1}{g},
	\]
	which is precisely
	\[
	\log\delta_g \asymp \frac{1}{g}.
	\]
	Equivalently, the minimal entropy is of order \(1/g\), and exponentiating
	the upper bound gives
	\[
	\delta_g \leq 11^{1/g}.
	\]
	In particular, \(\delta_g\to 1\) as \(g\to\infty\).
	
	Penner's upper bound was later improved by Bauer~\cite{bauer1992example},
	who showed that
	\[
	\log\delta_g \leq \frac{\log 6}{g}.
	\]
	Hironaka--Kin~\cite{hironaka2006small} constructed families with still
	smaller asymptotic constant, giving
	\[
	\log\delta_g \leq \frac{\log(2+\sqrt 3)}{g}.
	\]
	Their construction uses families of pseudo-Anosov braids and related
	covering constructions, producing explicit sequences of mapping classes
	whose entropies realize this improved upper bound.
	
	There are also exact results in small genus. For genus \(2\),
	Cho--Ham~\cite{cho2008minimum} proved that \(\delta_2\) is the largest
	real root of
	\[
	x^4-x^3-x^2-x+1,
	\]
	so numerically \(\delta_2\approx 1.72208\). Lanneau--Thiffeault
	\cite{lanneau2011minimum} studied the corresponding problem for small
	genera with orientable invariant foliations; in particular, they
	determined the minimum dilatations in genera \(3,4,\) and \(5\) in that
	orientable-foliation setting, and related the candidates to explicit
	reciprocal polynomials. For example, in genus \(3\) one obtains the largest
	real root of
	\[
	x^6-x^4-x^3-x^2+1,
	\]
	which is approximately \(1.40127\). These results show that, although
	Penner's theorem gives the correct asymptotic order, determining the exact
	value of \(\delta_g\) is a delicate arithmetic and geometric problem even
	for small genus.
	
	For punctured surfaces, related asymptotic questions have also been
	studied. Yazdi~\cite{yazdi2018lowerbound} proved lower bounds for
	pseudo-Anosov maps on \(S_{g,n}\), improving earlier quantitative
	estimates. He also addressed Penner's question on punctured surfaces:
	for fixed \(n\), the logarithm of the minimum dilatation on \(S_{g,n}\)
	again has order \(1/g\) as \(g\to\infty\)
	\cite{yazdi2019pseudo}.
			
	\subsection{Degrees of Stretch Factors and Multiplicity}
	
	Returning to the algebraic restrictions discussed in the previous
	sections, a natural question is which algebraic degrees can occur as
	degrees of pseudo-Anosov stretch factors on a fixed surface. Thurston
	observed that the degree of a pseudo-Anosov stretch factor on \(S_g\) is
	bounded above by \(6g-6\), the dimension of Teichmüller space. More
	precisely, the stretch factor appears as the Perron--Frobenius eigenvalue
	of a train-track transition matrix, whose relevant dimension is bounded by
	the dimension of the space of measured foliations.
	
	For pseudo-Anosov homeomorphisms with orientable invariant foliations, this
	bound improves to \(2g\). Indeed, in this case the invariant foliations
	determine cohomology classes which are eigenvectors of
	\[
	\phi^*:H^1(S_g;\mathbb R)\to H^1(S_g;\mathbb R)
	\]
	with eigenvalues \(\lambda\) and \(\lambda^{-1}\). Hence the minimal
	polynomial of \(\lambda\) divides the characteristic polynomial of
	\(\phi^*\), which has degree \(2g\). Therefore
	\[
	\deg(\lambda)\leq 2g.
	\]
	The sharpness and possible values of these degree bounds are studied in \cite{shin2014algebraic,strenner2017algebraic}. Strenner's results show that the degree restrictions are sharp: all degrees allowed by the homological and train-track bounds do occur. However, this is a statement only about possible degrees. It does not imply that every bi-Perron unit of one of those degrees is realized as a pseudo-Anosov dilatation.
	
	Indeed, the authors in~\cite{baik2019typical} demonstrate, at least for genus $g\geq10$, that  the proportion of bi-Perron numbers $\leq R$ of degree $\leq6g-6$ that actually arise as stretch factors for surface homeomorphisms on a surface of genus $\leq g$ limits to $0$ as $R \to \infty$. In other words, while stretch factors are always bi-Perron, a generic bi-Perron number of degree $6g-6$ is almost never a stretch factor on a surface of genus $g$, at least in high enough genera. 
	
	Another phenomenon, proved in \cite{baik2016constructing}, concerns
	multiplicity of stretch factors. For certain small genera
	\(g\in\{2,3,4,5\}\), there exist stretch factors which are realized by
	arbitrarily many distinct conjugacy classes of pseudo-Anosov mapping
	classes on \(S_g\). Equivalently, for every \(N\geq 1\), one can find a
	number \(\lambda>1\) and pairwise non-conjugate pseudo-Anosov mapping
	classes
	\[
	\phi_1,\ldots,\phi_N\in \operatorname{Mod}(S_g)
	\]
	such that
	\[
	\lambda(\phi_1)=\cdots=\lambda(\phi_N)=\lambda,
	\]
	giving arbitrarily large collections of distinct closed Teichmüller geodesics of equal length.

%Penner - Least Dil. in genus g.  $\delta_g \sim 1/11^g$\\
%
%Lanneau-Thiffeault - Found the actual $\delta_g$ in genus 2, 3, 4, 5, 6?. These were all Salem Numbers. \\
%
%Mehdi Yazdi - Hironaka-Kim - gave better bounds on $\delta_g$ than penner. \\
%
%Balazs Strenner - biPerron numbers of every even degree $2n \leq 2g$  occur on a surface of genus g. \\
%
%BRW - among all Biperron numbers less than R of degree less than 2g, the proportion of those which are stretch factors for surface homeos. on surface of genus $g$ or less with ``orientation'' limits to 0 as $R\to\infty$. Arbitrarily high multiplicity (equal length among geodesics) in Moduli space for genus 2,3,4 and 5.\\
	\section{Fibered 3-manifolds and Stretch Factors}
	\label{sec:fibered_3manifolds}

	The dilatation of a pseudo-Anosov mapping class is intrinsically a two-dimensional dynamical quantity, but its structural and asymptotic properties can be understood via three-dimensional manifolds. First, let us define:
	\begin{definition}
		Let $S = S_{g,n}$ be a connected, orientable surface of negative Euler characteristic, and let $\phi \in \mathrm{Mod}(S)$ be a pseudo-Anosov mapping class. The \emph{mapping torus}\index{mapping torus} of $\phi$ is the compact, orientable 3-manifold defined by the suspension space:
		\begin{equation}
			M_\phi = (S \times [0,1]) / \sim, \quad \text{where } (x, 1) \sim (\phi(x), 0).
		\end{equation}
	\end{definition}
	By Thurston’s foundational Hyperbolization Theorem for Fibered Manifolds (see for instance \cite{thurston1998notes}), the mapping torus $M_\phi$ admits a complete hyperbolic structure of finite volume if and only if the monodromy $\phi$ is pseudo-Anosov. 
	In this setup, the logarithm of the stretch factor, \(\log\lambda(\phi)\), which equals the topological entropy \(h(\phi)\), carries two geometric interpretations. On the one hand, within the mapping torus \(M_\phi\), the suspension of the stable and unstable measured foliations \((\mathcal F^\pm,\mu^\pm)\) gives rise to a pair of transverse two-dimensional singular foliations whose leaves are invariant under the suspension flow. On the other hand, in the language of Teichmüller theory, the mapping class $\phi$ acts as a translation along a unique invariant geodesic axis in the Teichmüller space $\mathcal{T}_{g,n}$ equipped with the Teichm\"{u}ller metric. The projection of this axis to the moduli space $\mathcal{M}_{g,n}$ forms a closed loop, and $\log \lambda(\phi)$ represents the precise Teichmüller length of this closed loop.
	
	Rather than fixing a topological surface $S$ and varying the mapping class $\phi$, this three-dimensional framework allows one to fix a single hyperbolic 3-manifold $M$ and analyze the infinite family of stretch factors generated by its various choices of fibrations over the circle.
	
	\subsection{The Thurston Norm and Fibered Cones}
	
	Let $M$ be a hyperbolic 3-manifold. To define the Thurston norm\index{norm!Thurston}, we let $H^1(M; \RR)$ be the first cohomology group of $M$, which is isomorphic via Poincaré duality to $H_2(M, \partial M; \RR)$. In \cite{thurston1986norm}, Thurston defined a semi-norm $\|\cdot\|_T$ on this space that measures the minimal topological complexity of surfaces dual to a given cohomology class.
	
	\begin{definition}[Thurston norm]
		For an integral cohomology class \(\alpha \in H^1(M; \ZZ)\), the \emph{Thurston semi-norm}\index{norm!Thurston} is defined by
		\begin{equation}
			\|\alpha\|_T = \min \left\{ \chi_-(S) \ \middle|\ \begin{array}{l} S \subset M \text{ is a properly embedded,} \\ \text{oriented surface dual to } \alpha \end{array} \right\},
		\end{equation}
		where $\chi_-(S)$ minimizes the complexity across all components $S_i$ of $S$ via:
		\begin{equation}
			\chi_-(S) = \sum_{i} \max\big(0, -\chi(S_i)\big).
		\end{equation}
		and \(S_i\) denotes the connected components of \(S\). The semi-norm extends uniquely to \(H^1(M;\RR)\) by homogeneity and continuity. It is
		usually referred to as the \emph{Thurston norm}.
	\end{definition}
	
	The unit norm ball, defined as 
	$$B_T = \{ \alpha \in H^1(M; \RR) \mid \|\alpha\|_T \leq 1 \},$$ 
	is a compact, convex, centrally symmetric polyhedron with finitely many faces and rational vertices. A top-dimensional face $F \subset \partial B_T$ is designated as a \emph{fibered face} if the open cone over $F$ contains a primitive integral class\footnote{An integral class is primitive if it cannot be written as $k\beta$ for any $\beta \in H^1(M; \mathbf{Z})$ and integer $k > 1$.} $\alpha \in H^1(M; \ZZ)$ that is represented by a globally defined fibration $\pi: M \to S^1$.	
	
	Thurston established a dichotomy for the top-dimensional faces of the norm ball: any such face is either entirely non-fibered or entirely fibered. Specifically, if a face $F$ contains at least one fibered class in its interior, then every primitive integral class lying within the open cone over the interior of $F$—defined as the \emph{fibered cone}\index{cone!fibered}
	\begin{equation}
		\mathcal{C}_F = \left\{ t \cdot x \in H^1(M; \RR) \mid t > 0, \, x \in \mathrm{int}(F) \right\}
	\end{equation}
	—corresponds to a fibration of $M$ over $S^1$. Importantly, for each such integral class $\alpha \in H^1(M; \ZZ)$ represented by a globally defined fibration $\pi \colon M \to S^1$, the corresponding fiber $S_\alpha = \pi^{-1}(e^{i\theta})$ is a connected, orientable surface that minimizes complexity in its homology class. This geometric minimality implies that $\|\alpha\|_T = -\chi(S_\alpha)$. Thus, a fibered hyperbolic $3$-manifold whose first betti number number is at least $2$ is guaranteed to produce a parametrized family of stretch factors arising from the primitive integral classes in the fibered cone.
	
%	Thurston established a powerful dichotomy: if a face $F$ contains one fibered class, then every primitive integral class lying within the open cone over $F$, defined as the \emph{fibered cone}\index{cone!fibered}:
%	\begin{equation}
%		\mathcal{C}_F = \left\{ t \cdot x \in H^1(M; \RR) \mid t > 0, \, x \in F \right\},
%	\end{equation}
%	corresponds to a unique fibration of $M$ over $S^1$. Crucially, for each such integral class $\alpha \in H^1(M; \mathbf{Z})$ that is represented by a globally defined fibration $\pi \colon M \to S^1$, the corresponding fiber \(S_\alpha=\pi^{-1}(e^{i\theta})\) is a connected, orientable surface that minimizes complexity in its homology class. This means $\|\alpha\|_T = -\chi(S_\alpha)$. The associated monodromy $\phi_\alpha: S_\alpha \to S_\alpha$ is guaranteed to be a pseudo-Anosov mapping class, thus providing an infinite, highly structured repository of stretch factors within a single manifold.
	
	\subsection{Fried's Theorem: Analytic Behavior of Stretch Factors}
	
	When the first Betti number satisfies $b_1(M) > 1$, a single fibered cone $\mathcal{C}_F$ indexes an infinite family of pseudo-Anosov maps on connected surfaces of arbitrarily high genus. To understand the global dynamics of the manifold, one wants to quantify how the stretch factor varies as a function of the chosen cohomology class $\alpha \in \mathcal{C}_F$. For each primitive integral class $\alpha \in \mathcal{C}_F \cap H^1(M; \ZZ)$, let $K(\alpha) = \lambda(\phi_\alpha) > 1$ denote the stretch factor of its associated surface monodromy. In \cite{fried1985growth}, Fried proved that the dynamical expansion behavior varies in a remarkably continuous and controlled fashion across the geometry of the cone.
	
	\begin{theorem}[Fried, 1985]
		\label{thm:fried_continuity}
		Let $M$ be a compact, oriented, fibered hyperbolic 3-manifold, and let $\mathcal{C}_F \subset H^1(M; \RR)$ be a fibered cone. The function $f: \mathcal{C}_F \cap H^1(M; \ZZ) \to \RR_{\geq 0}$ defined by 
		\begin{equation}
			f(\alpha) = \frac{1}{\log K(\alpha)}
		\end{equation}
		extends uniquely to a continuous, real-analytic, and strictly concave function on the entire open cone $\mathcal{C}_F$. Furthermore, $f$ is homogeneous of degree 1, satisfying:
		\begin{equation}
			f(t\alpha) = t f(\alpha) \quad \text{for all } t > 0,
		\end{equation}
		and fulfills the boundary condition:
		\begin{equation}
			\lim_{\alpha \to \partial \mathcal{C}_F} f(\alpha) = 0.
		\end{equation}
	\end{theorem}
	
	Fried's theorem imposes obstructions on the distribution of stretch factors. The homogeneity of degree 1 implies that the stretch factor scales along rays according to $\log K(t\alpha) = \frac{1}{t}\log K(\alpha)$. Furthermore, the strict concavity of \(f(\alpha)\) shows that minima of entropy on suitable compact affine sections of the cone occur in the interior. Because $f(\alpha) = 1/\log K(\alpha)$, this unique maximum for $f$ corresponds to the absolute global minimum for the entropy $\log K(\alpha)$. Conversely, as the class $\alpha$ approaches the topological boundary $\partial \mathcal{C}_F$ of the cone, $f(\alpha) \to 0$, causing the stretch factors to blow up exponentially to infinity ($\log K(\alpha) \to \infty$). This provides a direct mechanism for constructing low-dilatation mapping classes: one samples integral classes situated near the barycenter of a large fibered face.
	
	\subsection{McMullen's Teichmüller Polynomial}
	
	In \cite{mcmullen2000polynomial}, McMullen introduced an algebraic invariant capable of simultaneously encoding the stretch factors for all fibrations associated to a given fibered cone. This is achieved via the ``Teichm\"{u}ller polynomial''\index{polynomial!Teichm\"{u}ller}, which serves as a dynamical analogue to the classical Alexander polynomial\index{polynomial!Alexander}.\footnote{For a compact orientable \(3\)-manifold \(M\) obtained as a knot complement in the three sphere, and with \(b_1(M)>0\), one obtains a multivariable Alexander polynomial from the maximal free abelian cover, viewed as a module over the group ring of \(H_1(M;\mathbb Z)/\mathrm{torsion}\). McMullen showed that this polynomial determines an Alexander norm\index{norm!Alexander} on \(H^1(M;\mathbb R)\), which gives a lower bound for the Thurston norm\index{norm!Thurston} and agrees with it on fibered cones in many cases, see \cite{alexander1928topological,milnor1968infinite,mcmullen2002alexander}.}
	
	Let $G = H_1(M; \ZZ) / \mathrm{torsion} \cong \ZZ^b$, where $b = b_1(M)$ is the first Betti number of $M$. The group ring $\ZZ[G]$ is isomorphic to the ring of Laurent polynomials in $b$ variables, $\ZZ[t_1^{\pm 1}, \dots, t_b^{\pm 1}]$. 
	
	\begin{theorem}[McMullen, 2000]
		\label{thm:mcmullen_teich_poly}
		For each fibered face $F \subset H^1(M; \RR)$, there exists a polynomial $\Theta_F \in \ZZ[G]$, called the Teichmüller polynomial, that completely determines the stretch factors of all monodromies in the fibered cone $\mathcal{C}_F$. 
		
		Explicitly, any integral cohomology class $\alpha \in \mathcal{C}_F \cap H^1(M; \ZZ)$ specifies a group homomorphism $\alpha: G \to \ZZ$. This induces a ring homomorphism $\alpha_*: \ZZ[G] \to \ZZ[t^{\pm 1}]$ defined on the basis elements by $\alpha_*(g) = t^{ \alpha(g) }$. The stretch factor $K(\alpha)$ of the pseudo-Anosov monodromy is exactly equal to the largest real root of the single-variable specialized polynomial:
		\begin{equation}
			\alpha_*(\Theta_F)(t) = 0.
		\end{equation}
	\end{theorem}
	
	Algebraically, the Teichmüller polynomial computes the spectral radius of the transition matrices associated with the invariant train tracks of the monodromies. By encoding an infinite family of mapping class dilatations into a single multivariable polynomial, McMullen's framework shifts the study of least stretch factors from geometric topology into the study of the roots of specialized polynomials.
	
	\subsubsection{Return to the example} For the example in section \ref{example}, the mapping torus produces a hyperbolic $3$-manifold $M$ via Thurston's theorem. The homology class of the initial genus $2$ surface $S$ determines a fibered face of $M$. The homology group $G=\text{H}_1(M,\ZZ)/\text{torsion}$ is isomoprhic to $\ZZ^2$, where assume $(0,1)$ corresponds to the dual of the class $[S]$. The fibered cone is then determined by $\{(a,b)\in G:b>0, |a|<b\}$. The Teichm\"uller polynomial for this cone can be computed (up to multiplication by a unit in $\ZZ[t]$) to be \[\Theta=(u-1)(u^2-(t+5+\frac{1}{t})u+1).\] This means, for a point $(a,b)$ in the fibered cone, the pseudo-Anosov monodromy of the fibration corresponding to $(a,b)$ has stretch factor given by the largest root of the polynomial \[x^{2b}-(x^a+5+\frac{1}{x^a})x^b+1.\]
	
	\subsection{Fibered Faces and the Magic 3-Manifold}
	
	This 3-dimensional structural machinery provides the underlying mathematical foundation for understanding optimal upper and lower bounds on surfaces. A primary target for minimizing stretch factors is the \emph{magic 3-manifold}, denoted $N(1,1,1)$, which is the complement of the three-component chain link in $S^3$. The magic 3-manifold is highly symmetric, possesses a large Betti number ($b_1 = 3$), and contains several large fibered faces.
	
	The work of Kin--Takasawa~\cite{kinTakasawa2011magic} shows that the magic \(3\)-manifold organizes many small-entropy examples through its fibered faces: by studying the entropy function on a fibered face, they identify fibered classes whose monodromies realize the smallest known entropies for several families of punctured-sphere mapping classes. The work of Kin--Kojima--Takasawa~\cite{kinKojimaTakasawa2013minimal} extends this viewpoint to suitable Dehn fillings of the magic manifold, where certain fibered classes give monodromies on closed surfaces with very small dilatation.
	
	This should be understood as a fibered-face realization mechanism rather than as a statement that every small-dilatation example comes directly from the magic manifold itself. In particular, sequences of primitive integral classes in fibered cones, together with appropriate Dehn fillings, produce explicit fibrations whose monodromy maps are pseudo-Anosov homeomorphisms with controlled entropy. This perspective is closely related to the families constructed by Hironaka--Kin~\cite{hironaka2006small}, which give the closed-surface upper bound
	\[
		\log \lambda(\phi_g) \leq \frac{\log(2+\sqrt 3)}{g}.
	\]
	Equivalently, these constructions produce pseudo-Anosov homeomorphisms on closed genus-\(g\) surfaces whose entropies have the same \(1/g\) order as Penner's bound, but with a smaller explicit constant.

	For punctured surfaces, Yazdi~\cite{yazdi2019pseudo} verified Penner’s long-standing conjecture that for a fixed number of punctures $n$, the minimal stretch factor obeys an exact $\asymp 1/g$ asymptotic decay. The main technique used 2-dimensional constructions, utilizing bipartite graphs and carefully engineered invariant train tracks to explicitly bound the Perron--Frobenius eigenvalues of the resulting incidence matrices.

	\subsection{Explicit Sequences of Teichm\"uller Polynomials}
	
	Recent work by Rafiqi~\cite{rafiqi2025teichmullerpolynomial} provides an explicit formula for McMullen's Teich\-m\"ul\-ler polynomial in a structured setting: the mapping tori of pseudo-Anosov homeomorphisms whose invariant foliations are orientable, possess exactly one singularity, and have all singular rays fixed. By translating these maps into permutations of integers and systematically employing canonical dual bases within the zippered rectangle formalism, Rafiqi explicitly computes the lifts required for McMullen's construction. 
	
	As an application, for every genus $g \geq 2$ and each integer $p \geq 0$, one can construct a pseudo-Anosov homeomorphism $f_{g,p}$ on a closed surface of genus $g$ such that the first Betti number of its mapping torus $M_{g,p}$ is precisely $b_1(M_{g,p}) = g$. The Teichm\"uller polynomial associated with the corresponding fibered face is given by:
	\begin{equation}\label{eq:polys}
		\begin{split}
			\Theta_{g,p}(t_1, \dots, t_{g-1}, u) = {} & (u-1)^{2g-3} \\
			& \times \left( u^2 - \left( \sum_{i=1}^{g-1} \left(t_i + t_i^{-1}\right) + 2g + p + 1 \right) u + 1 \right),
		\end{split}
	\end{equation}
	where the set $\{t_1, \dots, t_{g-1}, u\}$ forms a basis for $H_1(M_{g,p};\ZZ)/\mathrm{torsion}$, and $u$ represents the cohomology class corresponding to the initial monodromy $f_{g,p}$. 
	
	By specializing this multivariable polynomial along carefully chosen integral rays within the fibered cone, one can extract a family of mapping classes. This leads to the following result:
	
	\begin{proposition}\label{prop:bi_perron_units}
		Let $m \geq 2$ and let $a_1, \dots, a_{m-1}$ be non-negative integers such that the vector
		\[
		v = (\underbrace{1,\dots,1}_{a_{m-1}}, \, \underbrace{2,\dots,2}_{a_{m-2}}, \, \dots, \, \underbrace{m-1,\dots,m-1}_{a_1}, \, m)
		\]
		is a primitive integral vector. Then, for any integer $a_m \geq 3 + 2(a_1 + \dots + a_{m-1})$, the largest real root of the reciprocal polynomial
		\begin{equation*}
			x^{2m} - a_1 x^{2m-1} - \dots - a_{m-1} x^{m+1} - a_m x^m - a_{m-1} x^{m-1} - \dots - a_1 x + 1 = 0
		\end{equation*}
		is a bi-Perron unit and corresponds precisely to the stretch factor of a pseudo-Anosov map on a connected surface.
	\end{proposition}
	
	This proposition is significant because it realizes a positive proportion of all bi-Perron units in any given degree as actual pseudo-Anosov stretch factors. The primitivity condition imposed on the vector $v$ is extremely mild; it is automatically satisfied, for instance, whenever $a_1 \neq 0$, $a_{m-1} \neq 0$, or when $m$ is a prime number. The proposition is established by demonstrating that the integer vectors $v$ determine primitive cohomology classes that reside strictly within the fibered cones of the polynomials $\Theta_{g,p}$ for appropriate choices of the parameters $g$ and $p$.

\bibliographystyle{alpha}

\begin{thebibliography}{{Sch}08}
	
	\bibitem[Ale28]{alexander1928topological}
	J.~W. Alexander.
	\newblock Topological invariants of knots and links.
	\newblock {\em Transactions of the American Mathematical Society},
	30(2):275--306, 1928.
	
	\bibitem[Bau92]{bauer1992example}
	M.~Bauer.
	\newblock An upper bound for the least dilatation.
	\newblock {\em Transactions of the American Mathematical Society},
	330(1):361--370, 1992.
	
	\bibitem[BH92]{bestvina1992train}
	Mladen Bestvina and Michael Handel.
	\newblock Train tracks and automorphisms of free groups.
	\newblock {\em Annals of Mathematics}, 135(1):1--51, 1992.
	
	\bibitem[BRW16]{baik2016constructing}
	Hyungryul Baik, Ahmad Rafiqi, and Chenxi Wu.
	\newblock Constructing pseudo-anosov maps with given dilatations.
	\newblock {\em Geometriae Dedicata}, 180(1):39--48, 2016.
	
	\bibitem[BRW19]{baik2019typical}
	Hyungryul Baik, Ahmad Rafiqi, and Chenxi Wu.
	\newblock Is a typical bi-{Perron} algebraic unit a pseudo-{Anosov} dilatation?
	\newblock {\em Ergodic Theory Dyn. Syst.}, 39(7):1745--1750, 2019.
	
	\bibitem[CB88]{CassonBleiler1988Automorphisms}
	Andrew~J. Casson and Steven~A. Bleiler.
	\newblock {\em Automorphisms of Surfaces after Nielsen and Thurston}, volume~9
	of {\em London Mathematical Society Student Texts}.
	\newblock Cambridge University Press, Cambridge, 1988.
	
	\bibitem[CH08]{cho2008minimum}
	Ji-Young Cho and Ji-Young Ham.
	\newblock The minimal dilatation of a genus-two surface.
	\newblock {\em Experimental Mathematics}, 17(3):257--267, 2008.
	
	\bibitem[Euc56]{euclidBookX}
	Euclid.
	\newblock {\em The Thirteen Books of Euclid's Elements}, volume~3, chapter Book
	X: {I}ncommensurable Magnitudes, pages 1--254.
	\newblock Dover Publications, New York, second edition, 1956.
	
	\bibitem[FLP12]{FathiLaudenbachPoenaru2012ThurstonWork}
	Albert Fathi, Fran{\c c}ois Laudenbach, and Valentin Po{\'e}naru.
	\newblock {\em Thurston's Work on Surfaces}, volume~48 of {\em Mathematical
		Notes}.
	\newblock Princeton University Press, Princeton, NJ, 2012.
	
	\bibitem[FM11]{farb/margalit/MCG}
	Benson {Farb} and Dan {Margalit}.
	\newblock {\em {A primer on mapping class groups.}}
	\newblock Princeton, NJ: Princeton University Press, 2011.
	
	\bibitem[Fri85]{fried1985growth}
	David Fried.
	\newblock Growth rate of surface homeomorphisms and flow equivalence.
	\newblock {\em Ergodic Theory and Dynamical Systems}, 5(4):539--563, 1985.
	
	\bibitem[Fro12]{frobenius/matrizen/aus/nicht/negativen/elementen}
	G.~Frobenius.
	\newblock {\"U}ber {Matrizen} aus nicht negativen {Elementen}.
	\newblock {\em Berl. Ber.}, 1912:456--477, 1912.
	
	\bibitem[HK06a]{hironaka2006small}
	Eriko Hironaka and Eiko Kin.
	\newblock A family of pseudo-{A}nosov braids with small dilatation.
	\newblock {\em Algebraic \& Geometric Topology}, 6:699--738, 2006.
	
	\bibitem[HK06b]{hironaka2006family}
	Eriko Hironaka and Eiko Kin.
	\newblock A family of pseudo-anosov braids with small dilatation.
	\newblock {\em Algebraic \& Geometric Topology}, 6(2):699--738, 2006.
	
	\bibitem[HS00]{hindry2000diophantine}
	Marc Hindry and Joseph~H. Silverman.
	\newblock {\em Diophantine Geometry: An Introduction}, volume 201 of {\em
		Graduate Texts in Mathematics}.
	\newblock Springer-Verlag, New York, 2000.
	
	\bibitem[Ker83]{Kerckhoff1983NielsenRealization}
	Steven~P. Kerckhoff.
	\newblock The nielsen realization problem.
	\newblock {\em Annals of Mathematics}, 117(2):235--265, 1983.
	
	\bibitem[KKT13]{kinKojimaTakasawa2013minimal}
	Eiko Kin, Sadayoshi Kojima, and Mitsuhiko Takasawa.
	\newblock Minimal dilatations of pseudo-{A}nosovs generated by the magic
	3-manifold and their asymptotic behavior.
	\newblock {\em Algebraic \& Geometric Topology}, 13(6):3537--3602, 2013.
	
	\bibitem[KT11]{kinTakasawa2011magic}
	Eiko Kin and Mitsuhiko Takasawa.
	\newblock Pseudo-{A}nosov braids with small entropy and the magic 3-manifold.
	\newblock {\em Communications in Analysis and Geometry}, 19(4):705--758, 2011.
	
	\bibitem[Lin84]{lind/entropies/of/markov/shifts}
	D.~A. Lind.
	\newblock The entropies of topological {M}arkov shifts and a related class of
	algebraic integers.
	\newblock {\em Ergodic Theory Dynam. Systems}, 4(2):283--300, 1984.
	
	\bibitem[LT11]{lanneau2011minimum}
	Erwan Lanneau and Jean-Luc Thiffeault.
	\newblock On the minimum dilatation of pseudo-{A}nosov homeomorphisms on
	surfaces of small genus.
	\newblock {\em Annales de l'Institut Fourier}, 61(1):105--144, 2011.
	
	\bibitem[McM00]{mcmullen2000polynomial}
	Curtis~T. McMullen.
	\newblock Polynomial invariants for fibered 3-manifolds and teichm{\"u}ller
	geodesics for foliations.
	\newblock {\em Annales Scientifiques de l'{\'E}cole Normale Sup{\'e}rieure},
	33(4):519--560, 2000.
	
	\bibitem[McM02]{mcmullen2002alexander}
	Curtis~T. McMullen.
	\newblock The {Alexander} polynomial of a \(3\)-manifold and the {Thurston}
	norm on cohomology.
	\newblock {\em Annales scientifiques de l'{\'E}cole Normale Sup{\'e}rieure},
	35(2):153--171, 2002.
	
	\bibitem[Mil68]{milnor1968infinite}
	John Milnor.
	\newblock Infinite cyclic coverings.
	\newblock In J.~G. Hocking, editor, {\em Conference on the Topology of
		Manifolds}, pages 115--133. Prindle, Weber \& Schmidt, Boston, 1968.
	
	\bibitem[MS80]{misiurewicz/entropy/of/piecewise/linear/maps}
	M.~Misiurewicz and W.~Szlenk.
	\newblock Entropy of piecewise monotone mappings.
	\newblock {\em Studia Math.}, 67(1):45--63, 1980.
	
	\bibitem[MT88]{milnor/thurston/iterated/maps/of/the/interval}
	John Milnor and William Thurston.
	\newblock On iterated maps of the interval.
	\newblock In {\em Dynamical systems ({C}ollege {P}ark, {MD}, 1986--87)}, volume
	1342 of {\em Lecture Notes in Math.}, pages 465--563. Springer, Berlin, 1988.
	
	\bibitem[Neu99]{neukirch/ant}
	J{\"u}rgen Neukirch.
	\newblock {\em Algebraic number theory}, volume 322 of {\em Grundlehren der
		Mathematischen Wissenschaften [Fundamental Principles of Mathematical
		Sciences]}.
	\newblock Springer-Verlag, Berlin, 1999.
	\newblock Translated from the 1992 German original and with a note by Norbert
	Schappacher, With a foreword by G. Harder.
	
	\bibitem[OP20]{OhshikaPapadopoulos2020GlimpseThurston}
	Ken'ichi Ohshika and Athanase Papadopoulos.
	\newblock A glimpse into thurston's work.
	\newblock In Ken'ichi Ohshika and Athanase Papadopoulos, editors, {\em In the
		Tradition of Thurston: Geometry and Topology}, pages 1--58. Springer, Cham,
	2020.
	
	\bibitem[Pap24]{Papadopoulos2024FootnotesThurstonNotes}
	Athanase Papadopoulos.
	\newblock Some footnotes on thurston's notes the geometry and topology of
	3-manifolds.
	\newblock In Ken'ichi Ohshika and Athanase Papadopoulos, editors, {\em In the
		Tradition of Thurston III: Geometry and Dynamics}, pages 423--447. Springer,
	Cham, 2024.
	
	\bibitem[Pen91]{penner1991bounds}
	Robert~C. Penner.
	\newblock Bounds on least dilatations.
	\newblock {\em Proceedings of the American Mathematical Society},
	113(2):443--450, 1991.
	
	\bibitem[Per07]{perron/theorie/jacobischen/kettenbruchalgorithmus}
	Oskar Perron.
	\newblock Grundlagen f\"ur eine {T}heorie des {J}acobischen
	{K}ettenbruchalgorithmus.
	\newblock {\em Math. Ann.}, 64(1):1--76, 1907.
	
	\bibitem[PH92]{penner-train-tracks}
	R.~C. Penner and J.~L. Harer.
	\newblock {\em Combinatorics of train tracks}, volume 125 of {\em Annals of
		Mathematics Studies}.
	\newblock Princeton University Press, Princeton, NJ, 1992.
	
	\bibitem[Raf25]{rafiqi2025teichmullerpolynomial}
	Ahmad Rafiqi.
	\newblock A simple computation of {T}eichm{\"u}ller polynomials from integer
	permutations.
	\newblock {\em arXiv preprint arXiv:2505.06930}, 2025.
	
	\bibitem[HLW26]{chenxi2026}
	Hillen, P., Loving, M., Wu, C.
	\newblock {Every Weak Perron Number is an End-Periodic Stretch Factor}. 
	newblock {\em arXiv preprint arXiv:2603.20491}, 2026.
	
	\bibitem[RZ24]{rafiqi2024symplectic}
	Ahmad Rafiqi and Ayberk Zeytin.
	\newblock Invariant symplectic forms on number fields.
	\newblock {\em Turk. J. Math.}, 48(1):53--61, 2024.
	
	\bibitem[{Sch}08]{schoof/arakelov}
	Ren\'e {Schoof}.
	\newblock {Computing Arakelov class groups.}
	\newblock In {\em {Algorithmic number theory. Lattices, number fields, curves
			and cryptography}}, pages 447--495. Cambridge: Cambridge University Press,
	2008.
	
	\bibitem[Shi16]{shin2014algebraic}
	Hyunshik Shin.
	\newblock Algebraic degrees of stretch factors in mapping class groups.
	\newblock {\em Algebr. Geom. Topol.}, 16(3):1567--1584, 2016.
	
	\bibitem[Str17]{strenner2017algebraic}
	Bal{\'a}zs Strenner.
	\newblock Algebraic degrees of pseudo-{Anosov} stretch factors.
	\newblock {\em Geom. Funct. Anal.}, 27(6):1497--1539, 2017.
	
	\bibitem[Thu80]{thurston1998notes}
	William~P. Thurston.
	\newblock {\em The Geometry and Topology of Three-Manifolds}.
	\newblock Princeton University Course Notes, 1980.
	\newblock Electronic version 1.1 distributed by MSRI (2002).
	
	\bibitem[Thu86]{thurston1986norm}
	William~P. Thurston.
	\newblock A norm for the homology of 3-manifolds.
	\newblock {\em Memoirs of the American Mathematical Society}, 59(339):99--130,
	1986.
	
	\bibitem[Thu88]{Thurston1988GeometryDynamics}
	William~P. Thurston.
	\newblock On the geometry and dynamics of diffeomorphisms of surfaces.
	\newblock {\em Bulletin of the American Mathematical Society}, 19(2):417--431,
	1988.
	
	\bibitem[Thu14]{thurston/entropy/in/dimension/one}
	William~P. Thurston.
	\newblock Entropy in dimension one.
	\newblock In {\em Frontiers in complex dynamics}, volume~51 of {\em Princeton
		Math. Ser.}, pages 339--384. Princeton Univ. Press, Princeton, NJ, 2014.
	
	\bibitem[Thu23]{Thurston2023geometry/of/3/manifolds}
	William~P. Thurston.
	\newblock {\em The Geometry and Topology of Three-Manifolds}.
	\newblock American Mathematical Society, Providence, RI, 2023.
	\newblock With a preface by Steven P. Kerckhoff.
	
	\bibitem[Yaz18]{yazdi2018lowerbound}
	Mehdi Yazdi.
	\newblock Lower bound for dilatations.
	\newblock {\em J. Topol.}, 11(3):602--614, 2018.
	
	\bibitem[Yaz20]{yazdi2019pseudo}
	Mehdi Yazdi.
	\newblock Pseudo-{Anosov} maps with small stretch factors on punctured
	surfaces.
	\newblock {\em Algebr. Geom. Topol.}, 20(4):2095--2128, 2020.
	
	\bibitem[Yaz21]{yazdi2020lower}
	Mehdi Yazdi.
	\newblock Lower bound for the {P}erron--{F}robenious degrees of {P}erron
	numbers.
	\newblock {\em Ergodic Theory and Dynamical Systems}, 41(3):939--955, 2021.
	
	\bibitem[Yaz22]{yazdi2022nonnegative}
	Mehdi Yazdi.
	\newblock Non-negative integral matrices with given spectral radius and
	controlled dimension.
	\newblock {\em Ergodic Theory Dyn. Syst.}, 42(10):3246--3269, 2022.
	
\end{thebibliography}

\end{document}